\documentclass[sigconf]{acmart}
\usepackage{booktabs} 
\usepackage{amsmath}
\usepackage{algorithm}
\usepackage[noend]{algpseudocode}

\usepackage{caption}
\usepackage{graphicx}
\usepackage{subcaption}
\usepackage{float}
\usepackage{multirow}
\usepackage{array}
\usepackage[symbol]{footmisc}
\usepackage{lipsum}
\usepackage{multicol}
\usepackage{changepage}
\usepackage{adjustbox}

\graphicspath{ {images/} }



\setcopyright{rightsretained}

\acmDOI{10.1145/nnnnnnn.nnnnnnn} 
\acmISBN{978-x-xxxx-xxxx-x/YY/MM} 

\begin{document}
\title[MoHAEA]{Hybrid Adaptive Evolutionary Algorithm for \\ Multi-objective Optimization}

\author{Jeisson Prieto}
\orcid{1234-5678-9012}
\affiliation{%
  \institution{Universidad Nacional de Colombia \\ 
  Facultad de Ciencias\\
  Departamento de Matem\'aticas}
  \city{Bogot{\'a}} 
  \state{Colombia} 
}
\email{japrietov@unal.edu.co}

\author{Jonatan Gomez}
\affiliation{%
  \institution{Universidad Nacional de Colombia\\
  Facultad de Ingenier\'ia\\
  Dpto. Ing. Sistemas e Industrial}
  \city{Bogot{\'a}} 
  \country{Colombia}}
\email{jgomezpe@unal.edu.co}

\renewcommand{\shortauthors}{Prieto \& Gomez}

\begin{abstract}
The major difficulty in Multi-objective Optimization Evolutionary Algorithms (MOEAs) is how to find an appropriate solution that is able to converge towards the true Pareto Front with high diversity.  Most existing methodologies, which have demonstrated their niche on various practical problems involving two and three objectives, face significant challenges in the dependency of the selection of the EA parameters. Moreover, the process of setting such parameters is considered time-consuming, and several research works have tried to deal with this problem. This paper proposed a new Multi-objective Algorithm as an extension of the Hybrid Adaptive Evolutionary algorithm (HAEA) called MoHAEA. MoHAEA allows dynamic adaptation of the application of operator probabilities (rates) to evolve with the solution of the multi-objective problems combining the dominance- and decomposition-based approaches. MoHAEA is compared with four states of the art MOEAs, namely MOEA/D, pa$\lambda$-MOEA/D, MOEA/D-AWA, and NSGA-II on ten widely used multi-objective test problems. Experimental results indicate that MoHAEA outperforms the benchmark algorithms in terms of how it is able to find a well-covered and well-distributed set of points on the Pareto Front.

\end{abstract}

%
%



\maketitle

\section{Introduction}
Many real-world applications have more than one conflicting objectives to be simultaneously optimized. These problems are the so-called multi-objective optimization problems (MOPs). Multi-objective Evolutionary Algorithms (MOEAs) are shown to be suitable for solving various complex MOPs \cite{coello2006evolutionary}. Without any further information from a decision-maker, MOEAs algorithms are usually designed to meet two common but often conflicting goals: minimizing the distances between solutions and the Pareto Front (PF) (i.e., convergence) and maximizing the spread of solutions along the PF (i.e., diversity). 

However, there are some challenges brought by a large number of objectives \cite{li2014evolutionary}: (i) with the increasing number of objectives, almost all solutions in a population become nondominated with one another \cite{ishibuchi2008evolutionary}; (ii) with the increase of the objective space in size, the conflict between convergence and diversity becomes aggravated \cite{purshouse2007evolutionary}; and (iii) due to the computational efficiency consideration, the population size used in EMO algorithms cannot be arbitrarily large. To address this challenges there are two ways: (i) Pareto dominance-based approaches that uses the pareto dominance relation together with a neighbor density estimator to evaluate and approximate the Pareto Front \cite{zitzler2001spea2, deb2002fast, nebro2009smpso}; and (ii) Decomposition approaches, which decomposes a MOP into a set of subproblems and collaboratively optimizes them \cite{zhang2007moea, siwei2011multiobjective, qi2014moea}.

Although Evolutionary Algorithms (EA) have been used successfully in solving multi-objective optimization problems, the performance of this technique (measured in time consumed and solution quality) depends on the selection of the EA parameters \cite{eiben1999parameter, schaffer1989study}. A set of approaches, called Parameter Adaptation (PA) techniques tried to eliminate the parameter setting process by adapting parameters through the algorithm's execution \cite{davis1989adapting, tuson1998adapting, spears1995adapting}, as is the case of the Hybrid Adaptive Evolutionary Algorithm (HAEA) \cite{Gomez2004}. HAEA allows dynamic adaptation of the application of operator probabilities (rates) to evolve with the solution of the problem.

In this paper, a Multi-objective version of HAEA called MoHAEA is described and evaluated. MoHAEA combines the dominance- and decomposition-based approaches and uses the PA technique of HAEA, to tackle the three challenges posed in multi-objective optimization and to avoid the parameter setting process. The performance of MoHAEA is compared against some peer algorithms: MOEA/D \cite{zhang2007moea}, pa$\lambda$-MOEA/D \cite{siwei2011multiobjective}, MOEA/D-AWA \cite{qi2014moea}, and NSGA-II \cite{deb2002fast} over the ZDT and DTLZ multi-objective problems \cite{zitzler2000comparison, deb2002scalable} using the benchmark proposed in \cite{qi2014moea}.

The remainder of this paper is organized as follows. Section \ref{sec:Background} introduces the background knowledge of this paper.  Section \ref{sec:MoHAEA} describes, in detail, the proposed MoHAEA. Section \ref{sec:Experiments} presents the compared algorithms, quality indicators, and algorithm settings used for performance comparison. Section \ref{sec:Results} provides the experimental results and discussion. Finally, Section \ref{sec:Conclusion} outlines  some conclusions and future work lines.

\section{Background}
\label{sec:Background}

\subsection{Multi-objective Optimization Problem (MOP)}
The multi-objective optimization problem (MOP) can be mathematically defined as follows

\begin{equation}
\label{eq:MOP}
\begin{matrix}
    \text{\textbf{min} } f(x) = ( f_1(x), f_2(x), \dots, f_m(x))^T \\
    \text{\textbf{subject to }}  x \in \Omega
\end{matrix}
\end{equation}

\noindent where $x = (x_1, x_2, \dots , x_n)^T$ is a $n$-dimensional decision variable vector from the decision space $\Omega$; $f : \Omega \to \Theta \subseteq \mathbb{R}^m$ consists a set of $m$ objective functions that map $x$ from $n$-dimensional decision space $\Omega$ to $m$-dimensional objective space $\Theta$.

\begin{definition}{}
Given two decision vectors $x, y \in \Omega$, $x$ is said to Pareto \textbf{dominate} $y$, denoted by $x \prec y$, iff $f_i(x) \leq f_i(y)$, for every $i \in \{1, 2, \dots ,m\}$, and $f_j(x) < f_j(y)$, for at least one index $j \in \{1, 2, \dots ,m\}$.
\end{definition}

\begin{definition}{}
A decision vector $x^{*} \in \Omega$ is \textbf{Pareto optimal} iff there is no $x \in \Omega$ such that $x \prec x^{*}$.
\end{definition}

\begin{definition}
The \textbf{Pareto set (PS)} is defined as
\begin{equation}
    PS = \{x \in \Omega | x \text{ is Pareto Optimal}\}
\end{equation}
\end{definition}

\begin{definition}
The \textbf{Pareto Front (PF)} is defined as
\begin{equation}
    PF = \{f(x) \in \mathbb{R}^m | x \in \text{ PS  }\}
\end{equation}
\end{definition}

Since objectives in (\ref{eq:MOP}) conflicted each other, no point in $\Omega$ simultaneously minimizes all the objectives. The best trade-offs among the objectives can be defined in terms of PF. Then, The goal of MOEAs is to move the nondominated objective vectors toward PF (convergence), and also to generate a proper distribution of these vectors over the PF (diversity).

\subsection{Decomposition methods}
The decomposition-based method decomposes MOP into a set of subproblems and collaboratively optimizes them. Note that the decomposition concept is so general that either aggregation functions or simpler MOPs \cite{liu2013decomposition} can be used to form subproblems. Among them, weighted sum, weighted Tchebycheff, and boundary intersection approaches are very popular \cite{miettinen2012nonlinear}. In decomposition-based approaches, the subproblems are assigned to multiple predefined reference points for the collaborative optimization purpose. 

\begin{figure}[!ht]
    \centering
    \begin{subfigure}[t]{0.3\textwidth}
        \includegraphics[width=\textwidth, trim={0 0.4cm 0 0},clip]{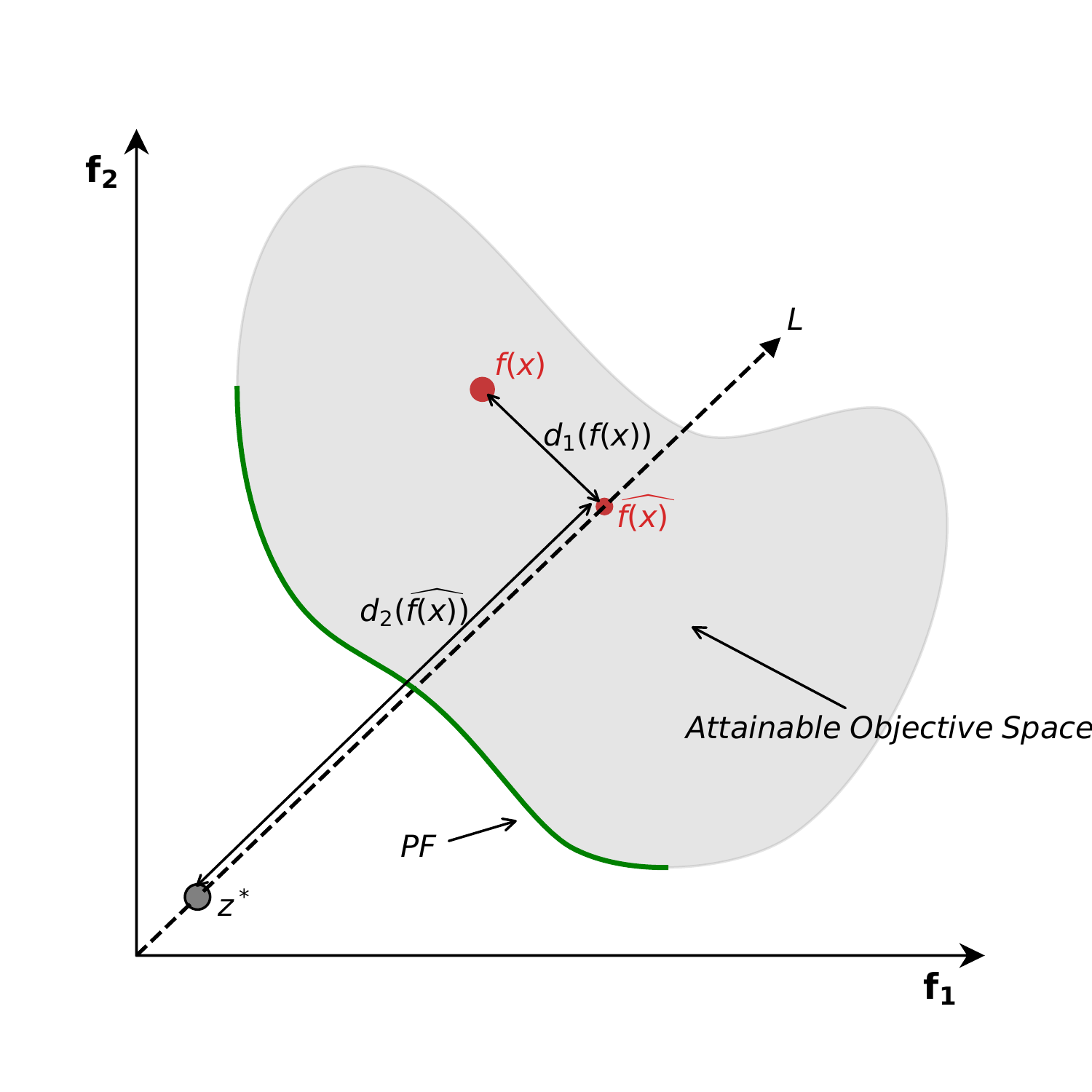}
    \end{subfigure}
    \caption{Penalty-based boundary intersection approach.}
    \label{fig:PBI}
\end{figure}

\begin{definition}
The ideal point $z^*$ or \textbf{reference point} is a vector $z^* =(z_1^*, z_2^*, \dots , z_m^*)^T$, with $z^* < \displaystyle \min_{x \in \Omega} f_i(x), i \in {1, \dots ,m}$
\end{definition}

One of the classic decomposition approach is the penalty-based boundary intersection (PBI) approach, which uses two distances: $d_1$ is used to evaluate the convergence of $x$ toward the PF and $d_2$ is a kind of measure for population diversity \cite{zhang2007moea}. Figure \ref{fig:PBI} presents a simple example to illustrate $d_1$ and $d_2$ of a solution $x$.

\begin{definition}
The penalty-based boundary intersection (\textbf{PBI}) approach is defined as
\begin{equation}
    \begin{matrix}
    \text{\textbf{min} } g(x| \lambda, z^*) = d_1 + \theta d_2 \\
    \text{\textbf{subject to }}  x \in \Omega
\end{matrix}
\end{equation}
\end{definition}

The proposed MoHAEA besides searching for non-dominated candidate solutions, also is trying to distribute on the approximated Pareto front by associating each one of those candidate solutions with a reference point. Since the above-created reference points are widely distributed on the entire normalized hyperplane, the obtained solutions are also likely to be widely distributed on or near the Pareto-optimal front.
    
\subsection{Hybrid Adaptive Evolutionary Algorithm}
Parameter adaptation techniques tried to eliminate the parameter setting process by adapting parameters through the algorithm's execution \cite{davis1989adapting, tuson1998adapting, spears1995adapting}. Parameter adaptation techniques can be roughly divided into centralized control techniques (central learning rule), decentralized control techniques, and hybrid control techniques \cite{eiben1999parameter}.

\begin{definition}
In the \textbf{centralized control techniques}, genetic operator rates are adapted according to a global learning rule based on the operator productivity through generations.
\end{definition}

\begin{definition}
In \textbf{decentralized control strategies}, genetic operator rates are encoded in the individual and are subject to the evolutionary process
\end{definition}

In \cite{Gomez2004}, a hybrid technique for parameter adaptation called Hybrid Adaptive Evolutionary Algorithm (HAEA) is proposed. In HAEA, each individual encodes its own operator rates and uses a randomized version of a learning rule mechanism for updating them is proposed. See Algorithm \ref{Al:HAEA}. 

In HAEA, each individual is "independently" evolved from the other individuals of the population, as in evolutionary strategies \cite{back1996evolutionary}. In each generation, every individual selects only one operator from the set of possible operators (line 8). Such an operator is selected according to the operator rates encoded into the individual. When a non-unary operator is applied, additional parents (the individual being evolved is considered a parent) are chosen according to any selection strategy, see line 9. Among the offspring produced by the genetic operator, only one individual is chosen as a child (line 11), and will take the place of its parent in the next population (line 17). In order to be able to preserve good individuals through evolution, HAEA compares the parent individual against the offspring generated by the operator (line 11). Therefore, an individual is preserved through evolution if it is better than all the possible individuals generated by applying the genetic operator.

\begin{algorithm}
\caption{HAEA \cite{Gomez2004}}
\label{Al:HAEA}
    \begin{flushleft}
    HAEA( $\lambda$, terminationCondition )
    \end{flushleft} 
    \begin{algorithmic}[1]
        \State $t_0 = 0$
        \State $\mathcal{P}_{0} = $ InitPopulation($\lambda$) 
        \While{( terminationCondition( $t, \mathcal{P}_{t}$ ) is \textit{false} )}
          \State $\mathcal{P}_{t+1} =\{ \}$
          \For{\textbf{each} $ind \in \mathcal{P}_{t}$}
          \State $rates$ = GetRates( $ind$ )
          \State $\delta$ = Random( $0,1$ ) \textit{ // learning rate}
          \State $operator$ = ChooseOperator( $operators, rates$ )
          \State $parents$ =  ParentSelection( $\mathcal{P}_{(t)}, ind$ ) 
          \State $offspring$ = ApplyOperator( $operator, parents$ )
          \State $child$ = Best( $offspring, ind$ )
          \If {($child_{fitness}$ $<$ $ind_{fitness}$)}
          \State $rates[oper]$ = ($1.0 + \delta$) $\times$  $rates[oper]$  \textit{ // reward}
          \Else
          \State $rates[oper]$ = ($1.0 - \delta$) $\times$  $rates[oper]$  \textit{ // punish}
          \EndIf
          \State NormalizeRates( $rates$ )
          \State SetRates( $child, rates$ )
          \State $\mathcal{P}_{t+1} = \mathcal{P}_{t+1} \cup \{child\}$
          \EndFor
            \State $t = t+1$
         \EndWhile
\end{algorithmic}
\end{algorithm}

\section{Proposed Algorithm: MoHAEA}
\label{sec:MoHAEA}

Algorithm \ref{Al:MoHAEA} presents the proposed Multi-Objective Hybrid Adaptive Evolutionary Algorithm (MoHAEA). This algorithm follows the main line of HAEA. 

\begin{algorithm}
\caption{MoHAEA}
\label{Al:MoHAEA}
    \begin{flushleft}
    MoHAEA($\lambda$, terminationCondition)
    \end{flushleft}  
    
    \begin{algorithmic}[1]
            \State $t_0 = 0$
            \State $\mathcal{P}_{0} = $ InitPopulation($\lambda$) 
                \While{( terminationCondition( $t, \mathcal{P}_{(t)}$ ) is \textit{false} )}
                  \State $\mathcal{P}_{t+1} =\{ \}$
                  \State UpdateDominanceCount($\mathcal{P}_t$)
                   \For{\textbf{each} $ind \in \mathcal{P}_{t}$}
                  \State $rates$ = GetRates( $ind$ )
                  \State $\delta$ = Random( $0,1$ ) \textit{ // learning rate}
                  \State $operator$ = ChooseOperator( $operators, rates$ )
                  \State $parents$ =  ParentSelection( $\mathcal{P}_{t}, ind$ ) 
                  \State $offspring$ = ApplyOperator( $operator, parents$ )
                  \State $child, rates$ = \textit{Best$^*$}( $offspring, ind, rates, oper$ )
                  \State NormalizeRates( $rates$ )
                  \State SetRates( $child, rates$ )
                  \State $\mathcal{P}_{t+1} = \mathcal{P}_{t} \cup \{child\}$
                  \EndFor
                    \State $t = t+1$
                 \EndWhile
    \end{algorithmic}
    
    \hphantom{--}
    \begin{flushleft}
    UpdateDominanceCount( $\mathcal{P}_{t}$ )
    \end{flushleft} 
    \begin{algorithmic}[1]
           \For{\textbf{each} $ind \in \mathcal{P}_{t}$}
                \State $ind_{\prec} = 0$ 
                \For{\textbf{each} $x \in \mathcal{P}_{t}$}
                    \If {$x \prec ind$}
                        \State $ind_{\prec} = ind_{\prec} + 1$
                    \EndIf
                \EndFor
          \EndFor
    \end{algorithmic}
    
    \hphantom{--}
    \begin{flushleft}
    \textit{Best*}( $offspring, ind, rates, oper$ ) 
    \end{flushleft}  
    \begin{algorithmic}[1]
        \State $x = KNN \left( ind , offspring \right)$ 
        \State UpdateFitness($x, ind$)
        \If {($x_{fitness} > ind_{fitness}$)}
            \State $x = ind$ \textit{// The offspring is not better to the parent}
        \EndIf    
        \If {($x_{fitness} < ind_{fitness}$)}
            \State $rates[oper]$ = ($1.0 + \delta$) $\times$  $rates[oper]$  \textit{ // reward}
          \Else
            \State $rates[oper]$ = ($1.0 - \delta$) $\times$  $rates[oper]$  \textit{ // punish}
        \EndIf    
        \State $x_{z^*} = ind_{z^*}$ \textit{// Reference point assignation}
        \State \textbf{return} $x, rates$
    \end{algorithmic}

    \hphantom{--}
    \begin{flushleft}
    UpdateFitness( $x, ind$ )
    \end{flushleft}  
    \begin{algorithmic}[1]
            \State $ind_{fitness} = d_1(f(ind))$
            \State $x_{fitness} = d_1(f(x))$
            \If {($x \prec ind$) \& ($\neg(ind \prec x$))}
                \State $ind_{fitness} = ind_{fitness} + 1$ \textit{// $ind$ non-dominance penalty}
            \EndIf
            \If {($ind \prec x$) \& ($\neg(x \prec ind$))}
                \State $x_{fitness} = x_{fitness} + 1$ \textit{// $x$ non-dominance penalty}
            \EndIf
    \end{algorithmic}
\end{algorithm}

\subsection{Individual encoding}
Following the encoding of the decentralized control adaptation techniques, the genetic operator rates are encoded into the individual with values following a uniform distribution $U[0, 1]$. Also, each individual encodes its reference point (randomly initialized using the \textit{InitPopulation} method).

\begin{figure}[!htb]
\adjustbox{width=0.45\textwidth}{
\begin{tabular}{|c|c|c|c|c|c|c|}
\hline
\multicolumn{1}{|c|}{\textbf{SOLUTION}} &  \multicolumn{3}{c|}{\textbf{OPERATORS}} & \multicolumn{1}{l|}{\textbf{REFERENCE POINT}} \\ \hline
$x$    & $Oper_1$  & $\dots$ & $Oper_k$ & $z^*$                                \\ \hline
{[}0.27 , $\dots$ , 0.01{]}  & 0.33      & $\dots$ & 0.10     & {[}0.50, $\dots$, 0.50{]}                     \\ \hline
\end{tabular}
}
\caption{Encoding operator probabilities and reference point in the chromosome}
\label{fig:encoding}
\end{figure}

 \begin{figure*}[t]
    \centering
    \begin{subfigure}[t]{0.3\textwidth}
        \includegraphics[width=\textwidth, trim={0 0.4cm 0 0},clip]{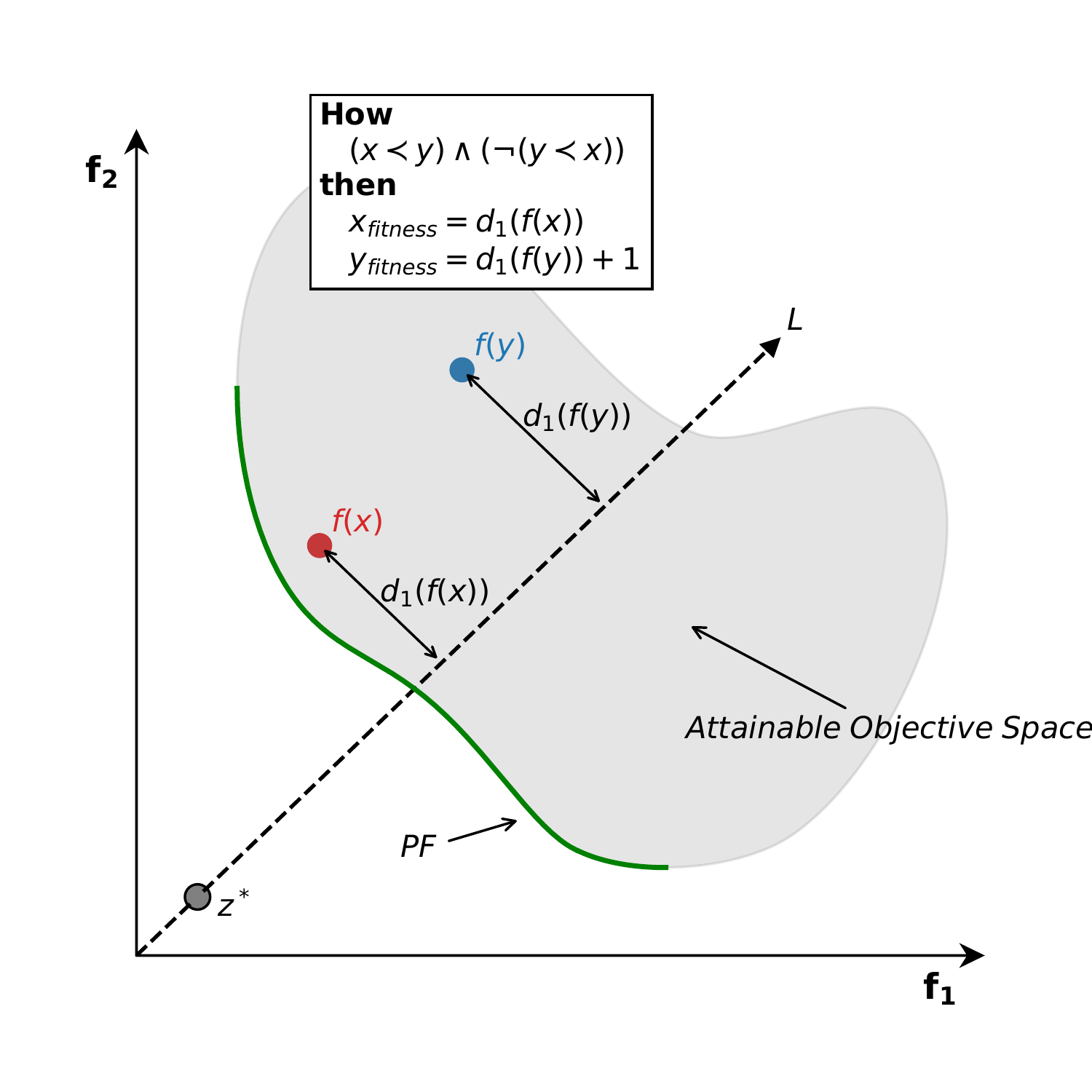}
        \caption{$(x \prec y) \land (\neg(y \prec x))$}
    \end{subfigure}
    \hphantom{---}
    \begin{subfigure}[t]{0.3\textwidth}
        \includegraphics[width=\textwidth, trim={0 0.4cm 0 0},clip]{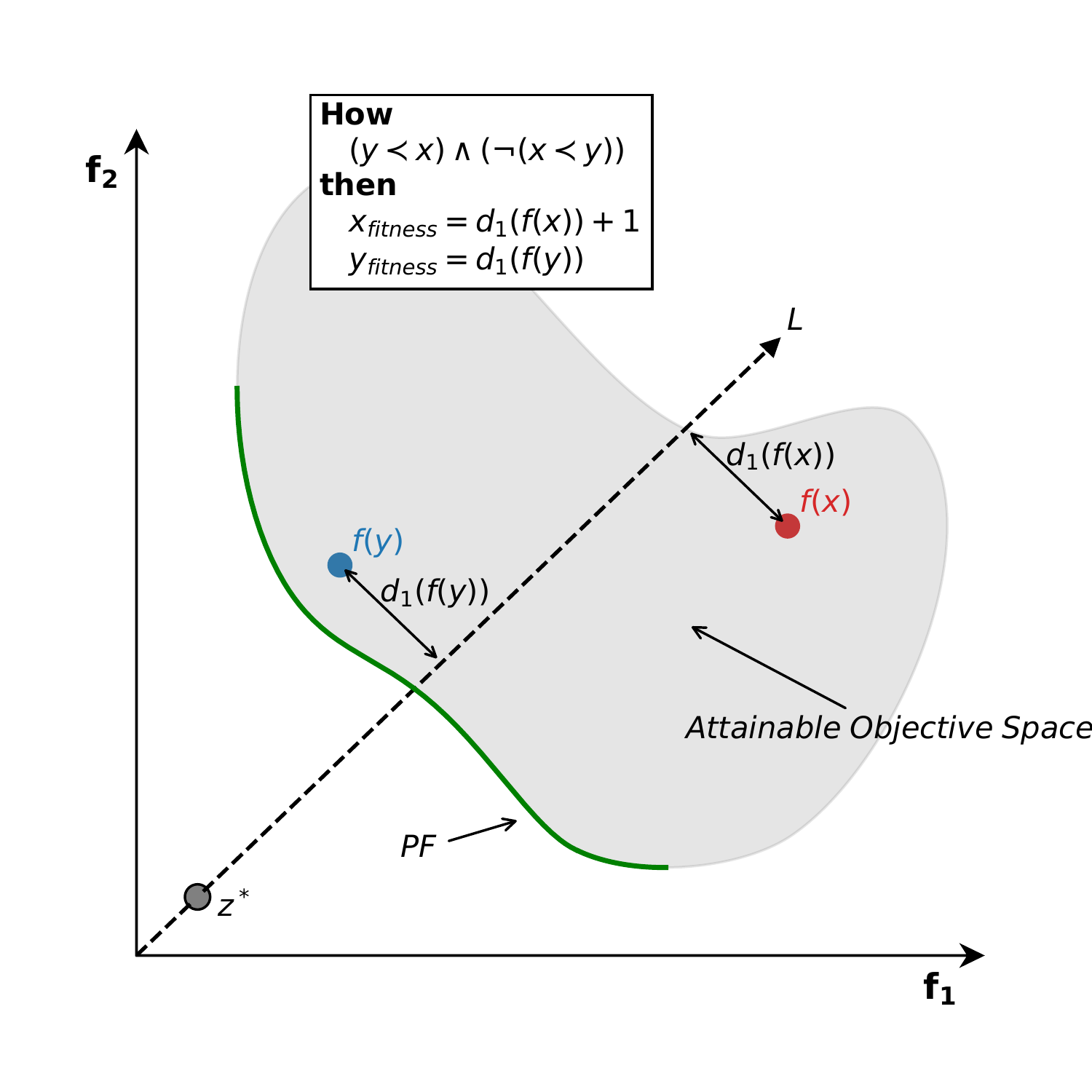}
        \caption{$(y \prec x) \land (\neg(x \prec y))$}
    \end{subfigure}   
    \hphantom{---}
    \begin{subfigure}[t]{0.3\textwidth}
        \includegraphics[width=\textwidth, trim={0 0.4cm 0 0},clip]{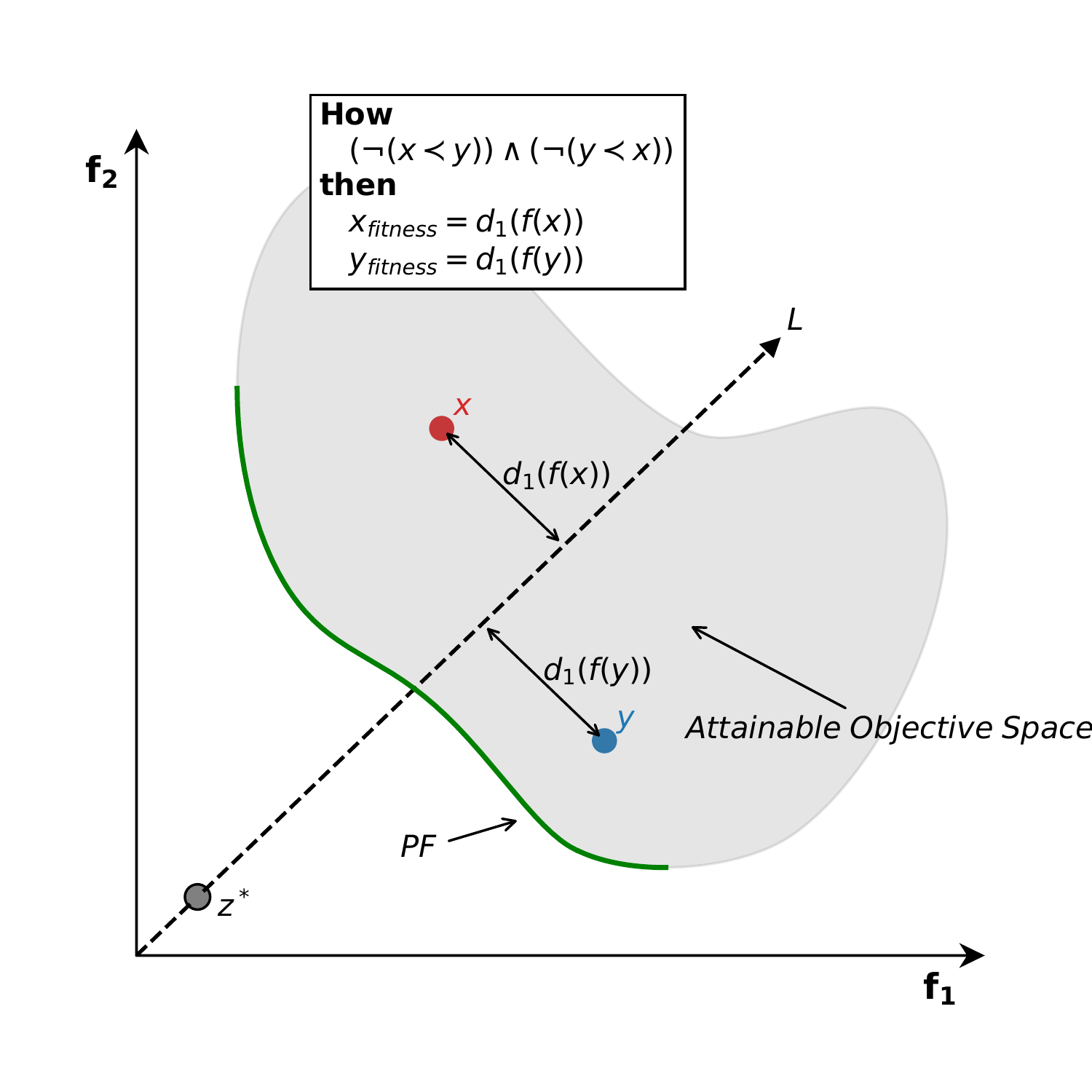}
        \caption{$(\neg(x \prec y)) \land (\neg(y \prec x))$}
    \end{subfigure}
    \caption{Proposed Dominance Penalty Approach for Fitness evaluation.}
    \label{fig:DominanceDecomposition}
\end{figure*}

\subsection{Selection Mechanism}
 In each generation, every individual selects only one operator from the set of possible operators using a roulette selection scheme (line 9). However, when a non-unary operator is applied, additional parents (the individual being evolved is considered a parent) are chosen according to a tournament selection of size 4 (line 10). This tournament is based on dominance; for this reason, it is necessary to update the dominance count for each individual every generation (line 4).

\subsection{Maintaining diversity}
For avoiding the crossover interaction problem, MoHAEA uses a modification of Deterministic Crowding (DC) proposed by Gomez \cite{gomez2004modal}. This modification allows MoHAEA to use the DC strategy not only with crossover and mutation but with different genetic operators. The algorithm determines the closest offspring to the individual using the \textit{K}-Nearest Neighbor ($KNN$) strategy \cite{altman1992introduction}, and it is used as a child in the generation (line 1, $Best^*$).

As indicated before, MoHAEA uses a predefined set of reference points to ensure diversity in obtained solutions. The Das and Dennis's systematic approach \cite{das1998normal} was used to choosing reference points. Then, each individual encodes its reference point $z^*$, and when an operator is applied, the individual passes its reference point to the offspring (line 9, $Best^*$).

\subsection{Dominance + Decomposition}
It follows the method the penalty-based boundary intersection (PBI) approach proposed in \cite{zhang2007moea}, which use two distances: $d_1$ is the distance between $f(x)$ and the reference line $L$ produced by $z^*$ and $d_2$ is the distance between and $z^*$ and $f(x)$. In MoHAEA, the distance $d_1$ is maintained (here, cosine distance), but the $d_2$ is now transformed into a "Dominance penalty", see Figure \ref{fig:DominanceDecomposition}. The fitness of the parent and the closest offspring are updated every generation using the dominance penalty approach (line 2, $Best^*$). 

\subsection{Computational Complexity of MoHAEA}
For each individual, MoHAEa requires a vector of ($op \times n \times m$) reals, where $op$ is the number of different genetic operators, $n$ is the number of decision variables, and $m$ is the number of objective functions (for saving the reference point). Thus, the space complexity of HAEA is linear concerning the number of objective functions (the population size of $N$ is considered a constant). Also, the time expended in calculating and normalizing the operator rates and fitness calculation is linear with respect at the population size $O(m \times N)$ (lines 12-14). The time complexity for updating of the Pareto Front for the selection strategy is $O(m \times N^2)$ (line 4).

The complexity of MoHAEA ($O(mN)$) in terms of fitness evaluations are lesser compared with NSGA-II ($O(m \times N^2)$), and the MOEA/D variants ($O(m \times N \times T)$ where $T$ is the neighborhood of the subproblem). Morover, MoHAEA reduces the extra parameters for the different compared algorithms such as the Neighborhood $T$, penalty $\theta$, or the $\sigma_{share}$.

\section{Experimental Setup}
\label{sec:Experiments}
In order to study the effectiveness of the proposed method (MoHAEA) a comparison with MOEA/D \cite{zhang2007moea}, pa$\lambda$-MOEA/D \cite{siwei2011multiobjective}, MOEA/D-AWA \cite{qi2014moea}, and NSGA-II \cite{deb2002fast} was made on different multi-objective test problems. Problems are a selection of five widely used bi-objective ZDT test instances \cite{zitzler2000comparison} and five tri-objective DTLZ problems \cite{deb2002scalable}.

\begin{table*}[]
\caption{Statistic IGD metric values of the Pareto-optimal solutions founded by the five compared algorithms on the ZDT and DTLZ problems. The numbers in parentheses are the computed deviations$^*$.}
\label{tb:comparison}
\begin{tabular}{lcccccc}
\hline
\multicolumn{1}{c}{Instance} & MoHAEA(SM) & MoHAEA(PM) & MOEA/D                                                            & pa$\lambda$-MOEA/D                                                        & MOEA/D-AWA                                                        &  NSGA-II                                                        \\\hline
ZDT1                         & \begin{tabular}[c]{@{}c@{}}3.779e-3\\ (2.208e-4)\end{tabular} &  \begin{tabular}[c]{@{}c@{}}3.833e-3\\ (2.678e-4)\end{tabular}      & \begin{tabular}[c]{@{}c@{}}4.739e-3\\ (3.973e-5)\end{tabular}     & \begin{tabular}[c]{@{}c@{}}\textbf{3.674e-3}\\ \textbf{(5.923e-5)}\end{tabular}     &
\begin{tabular}[c]{@{}c@{}}4.470e-3\\   (2.239e-4)\end{tabular}   &  \begin{tabular}[c]{@{}c@{}}4.696e-3\\ (1.435e-4)\end{tabular}  \\
ZDT2                         & \begin{tabular}[c]{@{}c@{}}\textbf{3.524e-3}\\ \textbf{(1.050e-4)}\end{tabular}  &\begin{tabular}[c]{@{}c@{}}3.661e-3\\ (1.647e-4)\end{tabular}      &  \begin{tabular}[c]{@{}c@{}}4.461e-3\\ (9.849e-5)\end{tabular}     & \begin{tabular}[c]{@{}c@{}}3.900e-3\\ (2.735e-4)\end{tabular}     & \begin{tabular}[c]{@{}c@{}}4.482e-3\\ (1.837e-4)\end{tabular}     & \begin{tabular}[c]{@{}c@{}}4.724e-3\\  (1.390e-4)\end{tabular} \\
ZDT3                         & \begin{tabular}[c]{@{}c@{}}3.817e-2\\ (1.106e-2)\end{tabular} &  \begin{tabular}[c]{@{}c@{}}3.596e-2\\ (8.065e-3)\end{tabular}     & \begin{tabular}[c]{@{}c@{}}1.362e-2\\  (1.574e-4)\end{tabular}    & \begin{tabular}[c]{@{}c@{}}9.737e-3\\  (6.636e-4)\end{tabular}    & \begin{tabular}[c]{@{}c@{}}6.703e-3\\  (4.538e-4)\end{tabular}    & \begin{tabular}[c]{@{}c@{}}\textbf{5.281e-3}\\ \textbf{(1.788e-4)}\end{tabular}  \\
ZDT4                         & \begin{tabular}[c]{@{}c@{}}\textbf{4.058e-3}\\ \textbf{(4.928e-4)}\end{tabular}  & \begin{tabular}[c]{@{}c@{}}1.102e-2\\ (2.152e-2)\end{tabular} &  \begin{tabular}[c]{@{}c@{}}4.692e-3\\     (1.339e-4)\end{tabular} & \begin{tabular}[c]{@{}c@{}}5.174e-3\\     (7.339e-4)\end{tabular} & \begin{tabular}[c]{@{}c@{}}4.238e-3\\     (3.102e-4)\end{tabular} & \begin{tabular}[c]{@{}c@{}}4.880e-3\\  (3.713e-4)\end{tabular} \\
ZDT6                         & \begin{tabular}[c]{@{}c@{}}6.707e-3\\ (1.145e-3)\end{tabular} & \begin{tabular}[c]{@{}c@{}}4.210e-2\\ (5.226e-3)\end{tabular}      & \begin{tabular}[c]{@{}c@{}}4.474e-3\\  (3.666e-4)\end{tabular}    & \begin{tabular}[c]{@{}c@{}}\textbf{3.601e-3}\\ \textbf{(4.250e-4)}\end{tabular}     & \begin{tabular}[c]{@{}c@{}}4.323e-3\\  (2.819e-4)\end{tabular}    &  \begin{tabular}[c]{@{}c@{}}4.261e-3\\ (2.255e-4)\end{tabular}  \\
DTLZ1                        & \begin{tabular}[c]{@{}c@{}}1.111e-2\\  (1.834e-3)\end{tabular} & \begin{tabular}[c]{@{}c@{}}\textbf{1.022e-2}\\  \textbf{(1.742e-3)}\end{tabular}  &  \begin{tabular}[c]{@{}c@{}}1.607e-2\\ (9.458e-4)\end{tabular}     & \begin{tabular}[c]{@{}c@{}}1.632e-2\\  (2.106e-3)\end{tabular}    & \begin{tabular}[c]{@{}c@{}}1.237e-2\\  (1.617e-3)\end{tabular}    & \begin{tabular}[c]{@{}c@{}}3.982e-2\\ (1.121e-3)\end{tabular}  \\
DTLZ2                        & \begin{tabular}[c]{@{}c@{}}\textbf{4.696e-3}\\  \textbf{(5.226e-4)}\end{tabular}  & \begin{tabular}[c]{@{}c@{}}5.284e-3\\  (6.971e-4)\end{tabular}      &  \begin{tabular}[c]{@{}c@{}}3.878e-2\\ (2.974e-4)\end{tabular}     & \begin{tabular}[c]{@{}c@{}}3.232e-2\\ (9.275e-4)\end{tabular}     & \begin{tabular}[c]{@{}c@{}}3.065e-2\\ (1.183e-4)\end{tabular}     & \begin{tabular}[c]{@{}c@{}}4.696e-2\\ (1.435e-3)\end{tabular}  \\
DTLZ3                        & \begin{tabular}[c]{@{}c@{}}1.036e-1\\  (9.888e-2)\end{tabular}   & \begin{tabular}[c]{@{}c@{}}7.389e-2\\  (2.221e-2)\end{tabular}      &  \begin{tabular}[c]{@{}c@{}}3.921e-2\\ (5.883e-4)\end{tabular}     & \begin{tabular}[c]{@{}c@{}}5.723e-2\\ (9.761e-2)\end{tabular}     & \begin{tabular}[c]{@{}c@{}}\textbf{3.196e-2}\\ \textbf{(8.036e-4)}\end{tabular}     & \begin{tabular}[c]{@{}c@{}}8.741e-2\\ (5.430e-2)\end{tabular}  \\
DTLZ4                        & \begin{tabular}[c]{@{}c@{}}1.960e-2\\  (1.367e-3)\end{tabular}  & \begin{tabular}[c]{@{}c@{}}\textbf{1.871e-2}\\  \textbf{(1.629e-3)}\end{tabular}      &  \begin{tabular}[c]{@{}c@{}}3.889e-2\\ (3.202e-4)\end{tabular}     & \begin{tabular}[c]{@{}c@{}}3.443e-2\\ (9.476e-3)\end{tabular}     & \begin{tabular}[c]{@{}c@{}}3.068e-2\\ (1.351e-4)\end{tabular}     & \begin{tabular}[c]{@{}c@{}}3.951e-2\\ (1.065e-3)\end{tabular}  \\
DTLZ6                        & \begin{tabular}[c]{@{}c@{}}1.808e-2\\  (1.599e-3)\end{tabular}  & \begin{tabular}[c]{@{}c@{}}\textbf{1.776e-2}\\  \textbf{(2.296e-3)}\end{tabular}      &  \begin{tabular}[c]{@{}c@{}}8.778e-2\\ (2.381e-3)\end{tabular}     & \begin{tabular}[c]{@{}c@{}}6.980e-2\\ (2.539e-3)\end{tabular}     & \begin{tabular}[c]{@{}c@{}}3.610e-2\\ (5.054e-3)\end{tabular}     &  \begin{tabular}[c]{@{}c@{}}4.156e-2\\ (1.483e-3)\end{tabular}  \\\hline
\end{tabular}

\small{$^*$ Values in bold indicate the best performing algorithm for the particular instance of a test problem.}
\end{table*}

\subsection{Experimental Setting}
In experimental studies, the comparison is carried on using the previous results of the algorithms in comparison: MOEA/D, pa$\lambda$-MOEA/D, MOEA/D-AWA, and NSGA-II reported in \cite{qi2014moea}. Then, the population size is set to $N = 100$ for the five bi-objective ZDT problems, and set to $N = 300$ for the given tri-objective DTLZ problems. Every algorithm stops when the number of function evaluations reaches the maximum number, 50,000 for the five bi-objective ZDT problems, and 75,000 for the five tri-objective DTLZ problems. 

As MOEA/D, pa$\lambda$-MOEA/D, MOEA/D-AWA, and NSGA-II \cite{qi2014moea}, MoHAEA uses the Simulated Binary Crossover (SBX) and Polynomial Mutation (PM). Distribution indexes for both SBX and PM are set to 20. Crossover rate is set to 1.00, while the mutation rate is set to $\frac{1}{n}$, where $n$ is the number of decision variables.

\renewcommand{\thefootnote}{\arabic{footnote}}
On the other hand, two different combinations of genetic operators were tested in MoHAEA: (i) Simulated Binary Crossover (SBX), Uniform \& Uniform (U\&U) operator, and Polynomial Mutation (PM) (MoHAEA(PM)); (ii) Simulated Binary Crossover (SBX), Uniform \& Uniform (U\&U) operator, and Shrink Mutation (SM) \cite{da2014simplex} (MoHAEA(SM)) \footnote{The distribution indexes in both SBX and PM are set to be 20. Uniform \& Uniform (U\&U) operator makes a Uniform Crossover operation generates two children, which are then modified by the Uniform mutation operator. SM with $\sigma = \frac{max-min}{20}$ , where max and min are the maximum and minimum values that can take each decision variable}. 

\subsection{Compared Algorithms}
To verify the proposed MoHAEA, the following four state-of-the-art algorithms for multi-objective functions are considered as peer algorithms.

\begin{itemize}
    \item \textbf{MOEA/D \cite{zhang2007moea}:} It is a representative of the decomposition-based method, the basic idea of MOEA/D is to decompose a MOP into a number of single-objective optimization subproblems through aggregation functions and simultaneously optimizes them.
    \item \textbf{ pa$\lambda$-MOEA/D \cite{siwei2011multiobjective}:} It is a based on decomposition with Pareto-adaptive weight vectors. This approach automatically adjusts the weight vectors by the geometrical characteristics of Pareto front. 
    \item \textbf{MOEA/D-AWA \cite{qi2014moea}}: An enhanced MOEA/D with weight vector initialization method and an elite population-based adaptive weight vector adjustment (AWA) strategy to address the MOPs with complex PFs.
    \item \textbf{NSGA-II \cite{deb2002fast}} The classic dominance approach strategy. The main feature is the use a fast nondominated sorting and crowded distance estimation procedure for comparing qualities of different solutions and selection.
\end{itemize}
 
\subsection{Performance Metric}
The Inverted Generational Distance ($IGD$) metric \cite{zitzler2003performance} was used to evaluate the performance of all compared algorithms in terms of convergence and diversity. Let $P^{*}$ be a set of evenly distributed points over the Pareto Front (in the objective space). Suppose that $P$ is an approximate set of the Pareto Front, the average distance from $P^{*}$ to $P$ is defined as:

\begin{equation}
    \label{eq:IGD}
    IGD(P^*, P) = \frac{\sum_{\nu \in P^*} d(\nu, P)}{|P^*|}
\end{equation}

\noindent where $d(\nu, P)$ is the minimum Euclidean distance between $\nu$ and the solutions in $P$. When $|P^{*}|$ is large enough, $IGD(P^{*}, P)$ can measure both the uniformity and the convergence of $P$. A low value of $IGD(P^{*}, P)$ indicates that $P$ is close to the Pareto Front and covers most of the whole Pareto Front.

For every benchmark algorithm, the $IGD$ metric is calculated using its population $\mathcal{P}$.

 \begin{figure*}[t]
    \centering
    \begin{subfigure}[t]{0.25\textwidth}
        \includegraphics[width=0.8\textwidth, trim={0 0.4cm 0 0},clip]{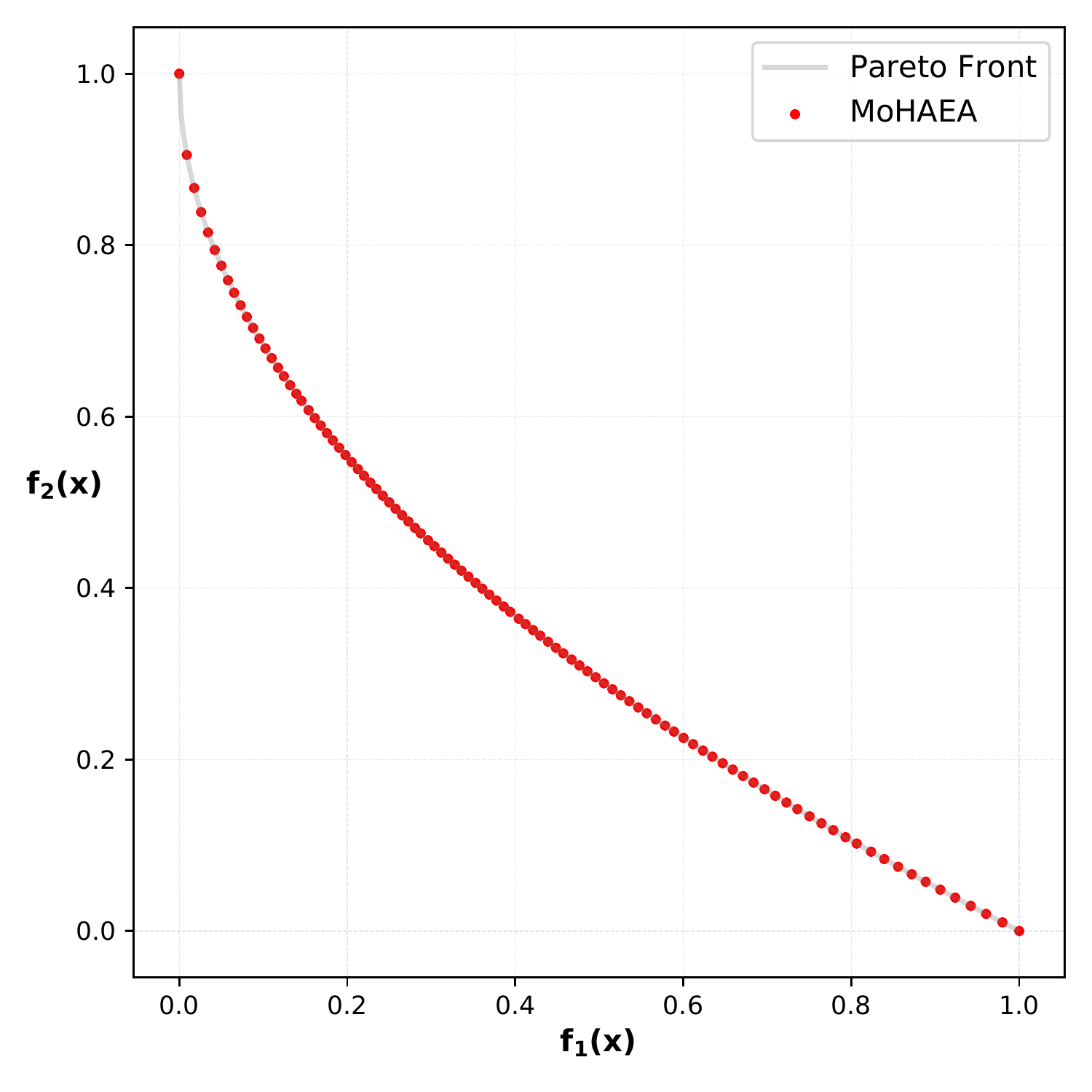}
        \caption{ZDT1}
    \end{subfigure}
    \hphantom{---}
    \begin{subfigure}[t]{0.25\textwidth}
        \includegraphics[width=0.8\textwidth, trim={0 0.4cm 0 0},clip]{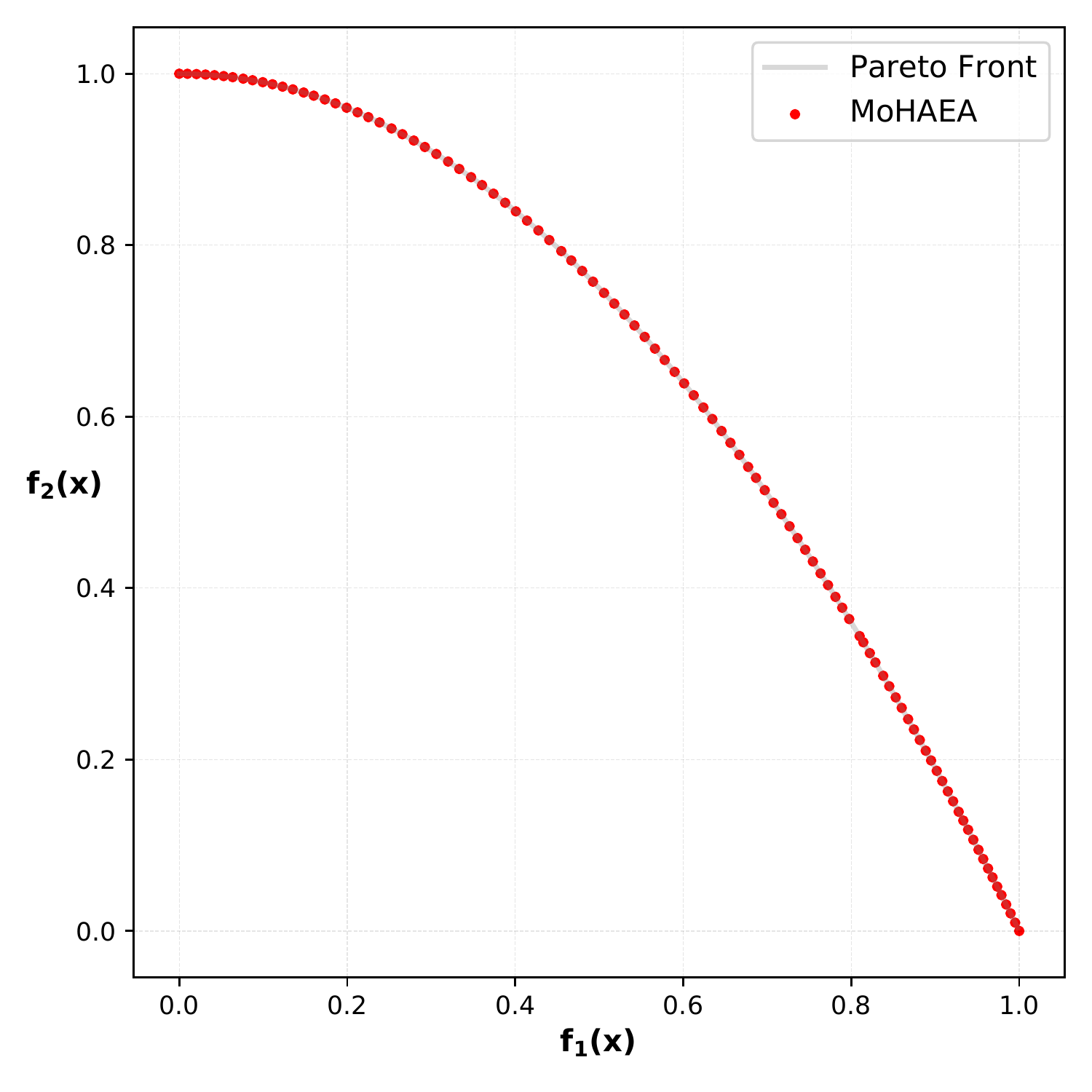}
        \caption{ZDT2}
    \end{subfigure}   
    \hphantom{---}
    \begin{subfigure}[t]{0.25\textwidth}
        \includegraphics[width=0.8\textwidth, trim={0 0.4cm 0 0},clip]{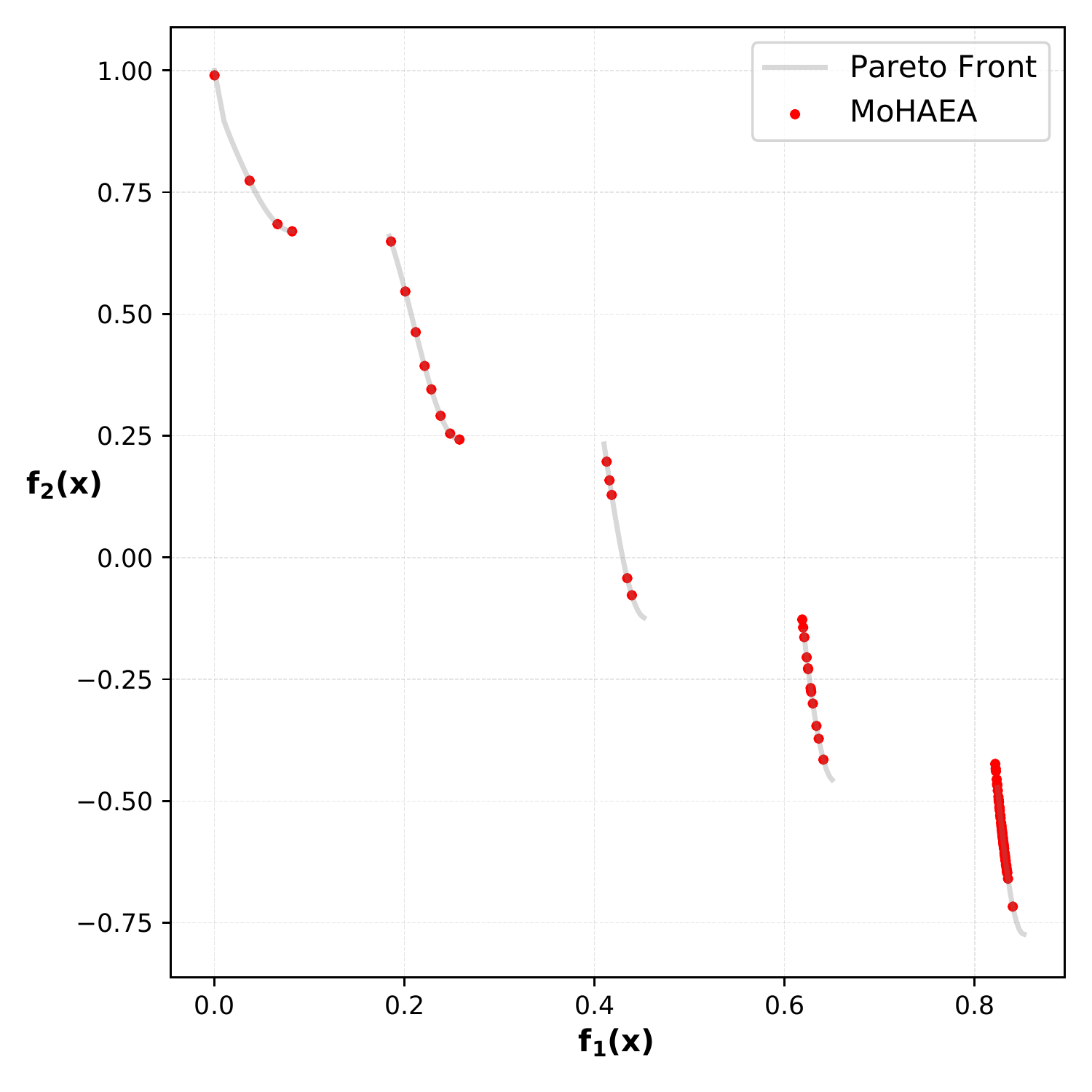}
        \caption{ZDT3}
    \end{subfigure}
    
    \begin{subfigure}[t]{0.25\textwidth}
        \includegraphics[width=0.8\textwidth, trim={0 0.4cm 0 0},clip]{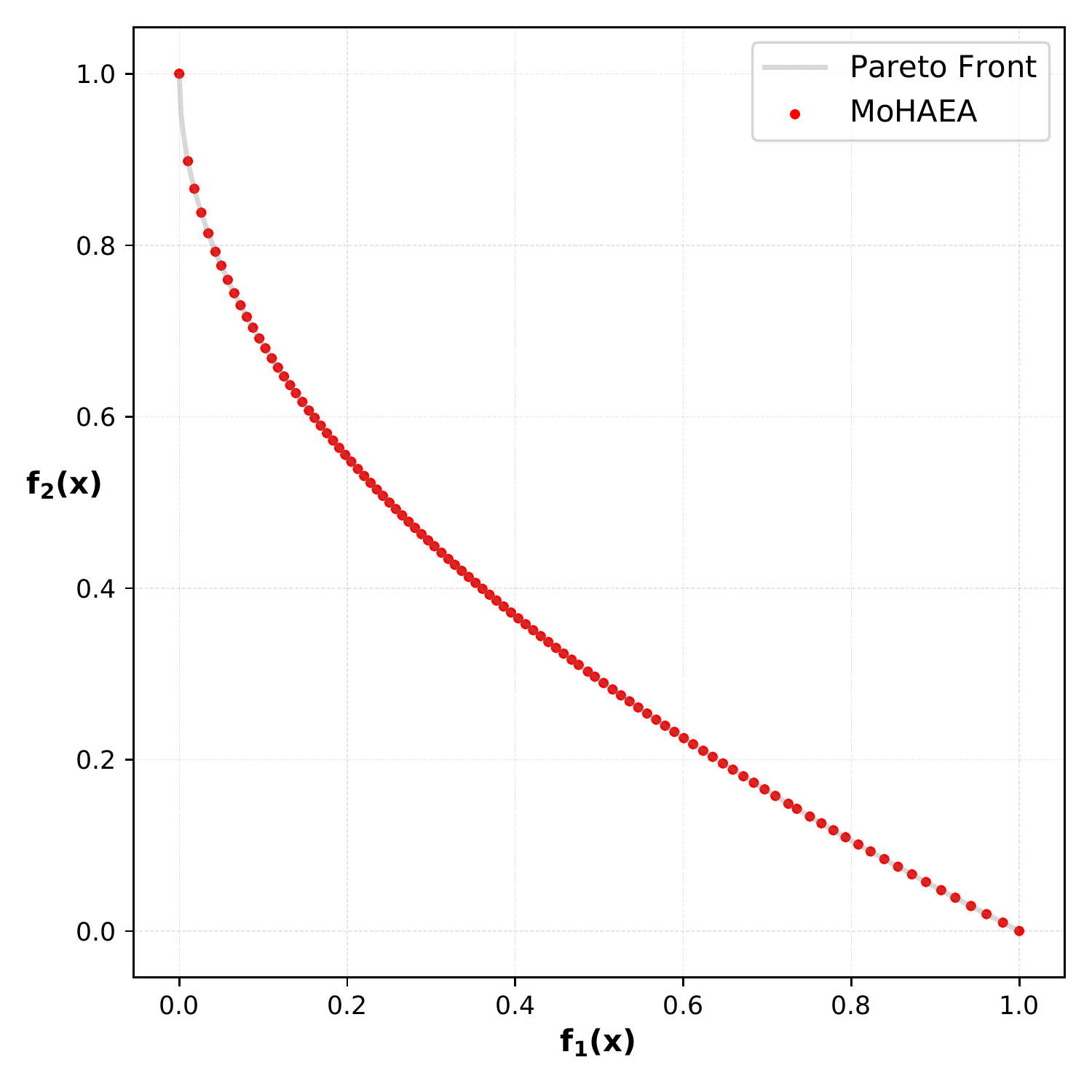}
        \caption{ZDT4}
    \end{subfigure}
    \begin{subfigure}[t]{0.25\textwidth}
        \includegraphics[width=0.8\textwidth, trim={0 0.4cm 0 0},clip]{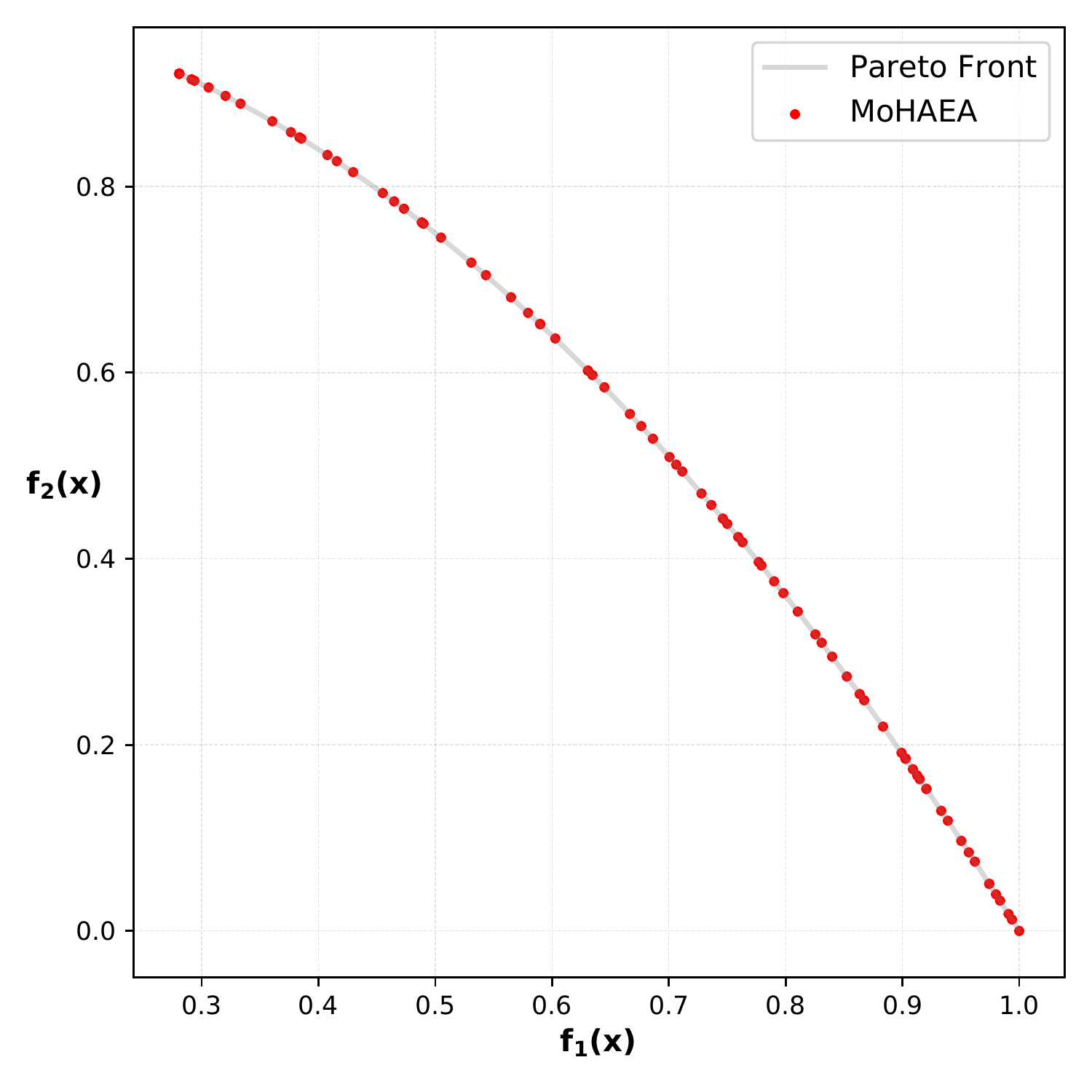}
        \caption{ZDT6}
    \end{subfigure}
    
    \begin{subfigure}[t]{0.25\textwidth}
        \includegraphics[width=\textwidth, trim={0.5cm 2.5cm 0 0.5cm},clip]{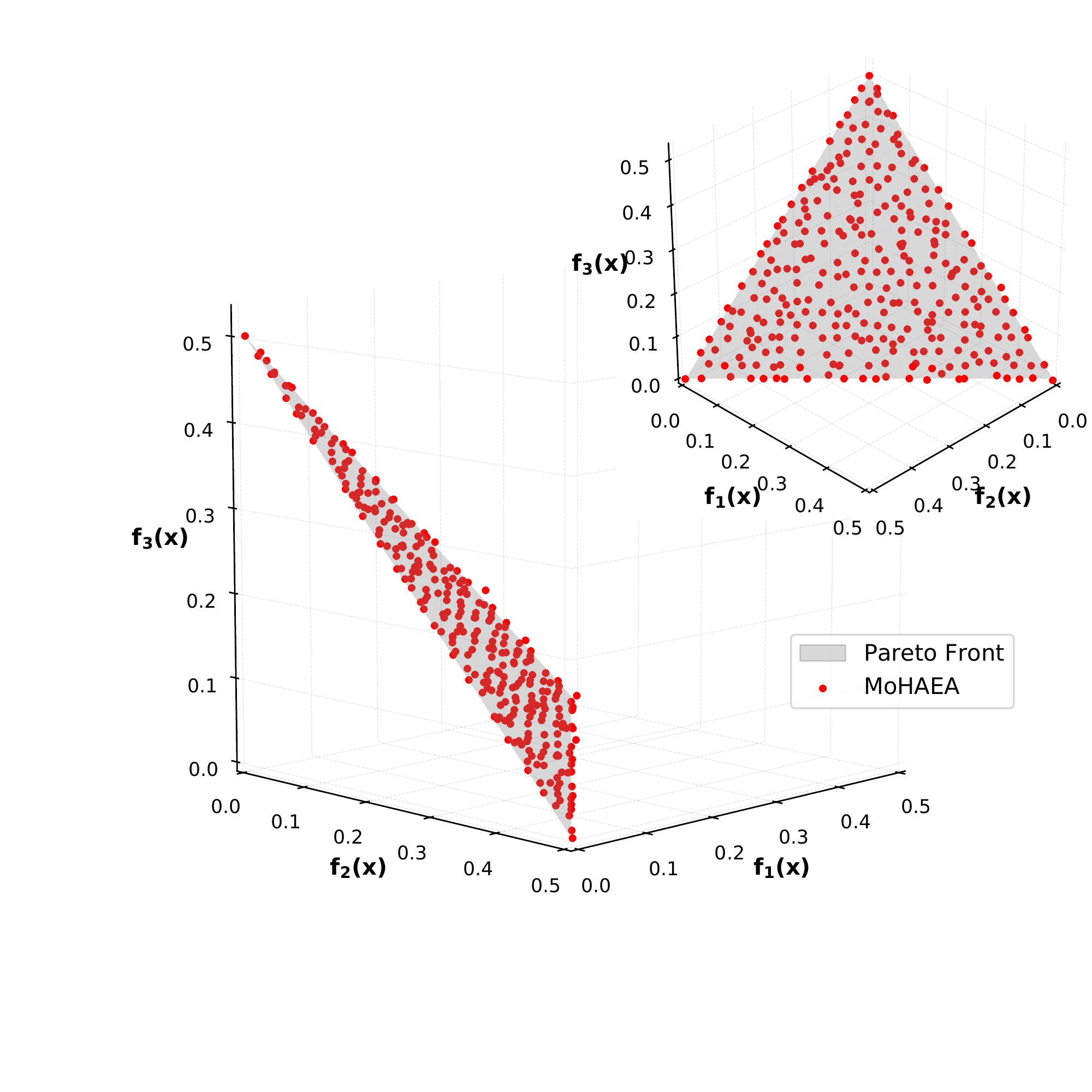}
        \caption{DTLZ1}
    \end{subfigure}
    \hphantom{---}
    \begin{subfigure}[t]{0.25\textwidth}
        \includegraphics[width=\textwidth, trim={0.5cm 2.5cm 0 0.5cm},clip]{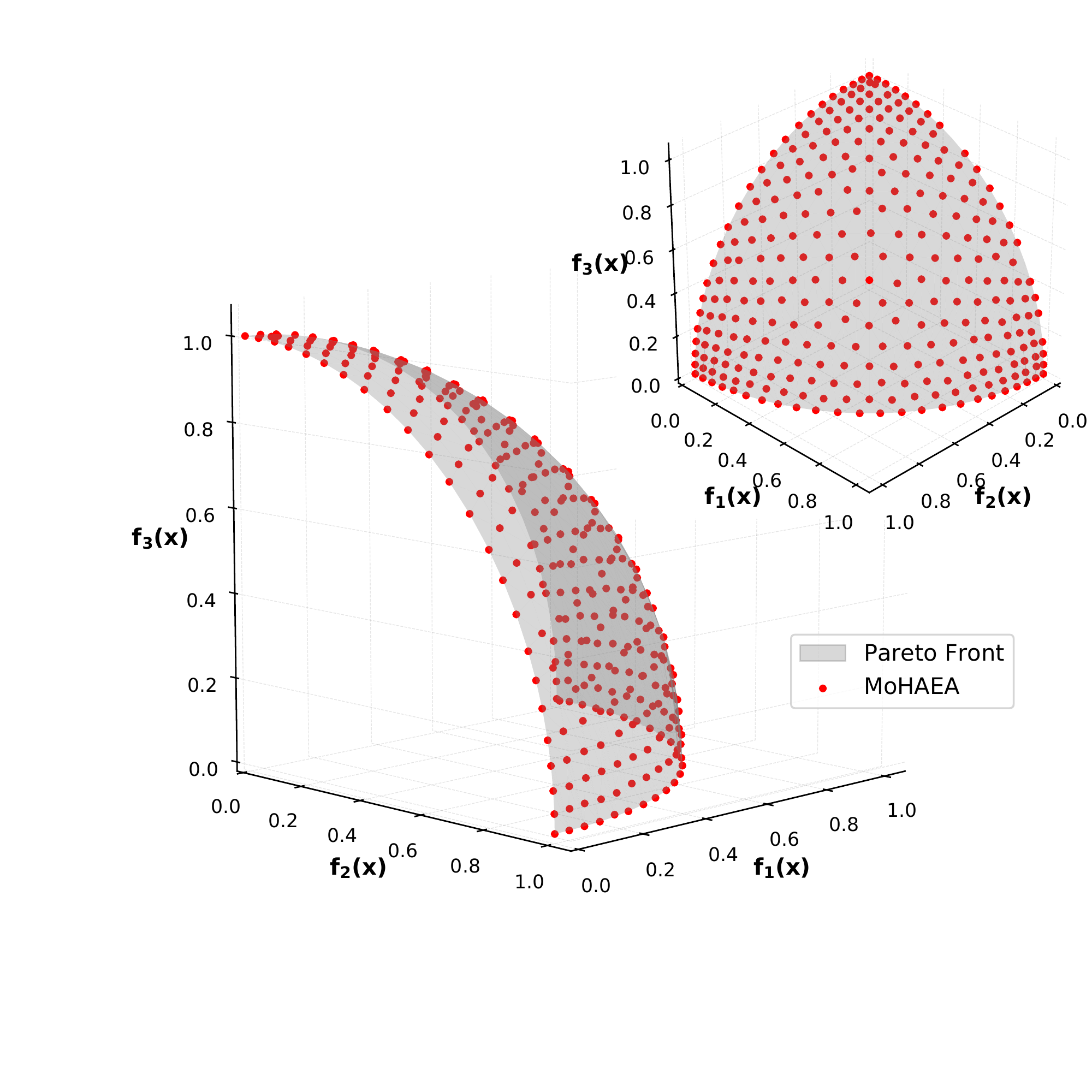}
        \caption{DTLZ2}
    \end{subfigure}   
    \hphantom{---}
    \begin{subfigure}[t]{0.25\textwidth}
        \includegraphics[width=\textwidth, trim={0.5cm 2.5cm 0 0.5cm},clip]{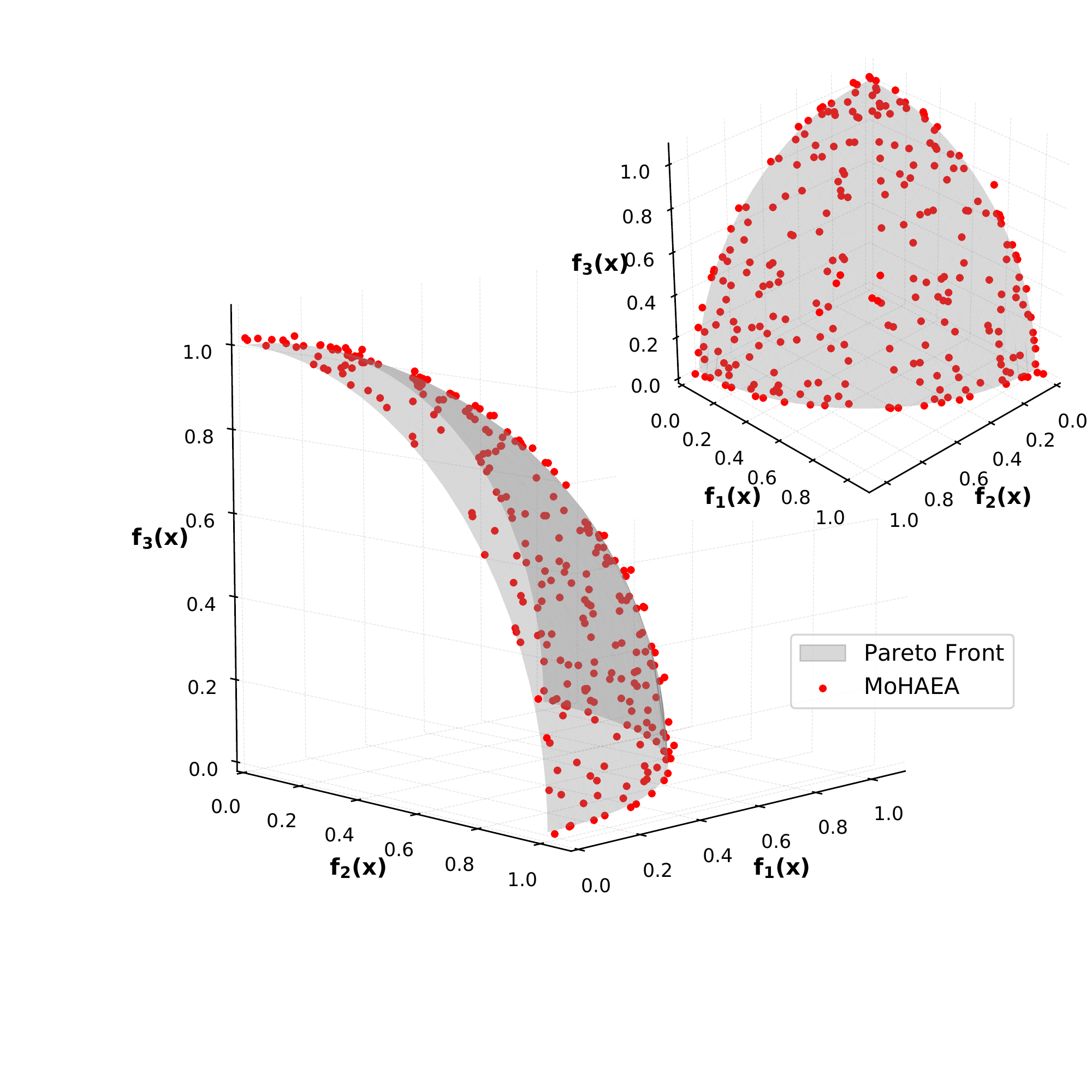}
        \caption{DTLZ3}
    \end{subfigure}
    
    \begin{subfigure}[t]{0.27\textwidth}
        \includegraphics[width=\textwidth, trim={0.5cm 2.5cm 0 0.5cm},clip]{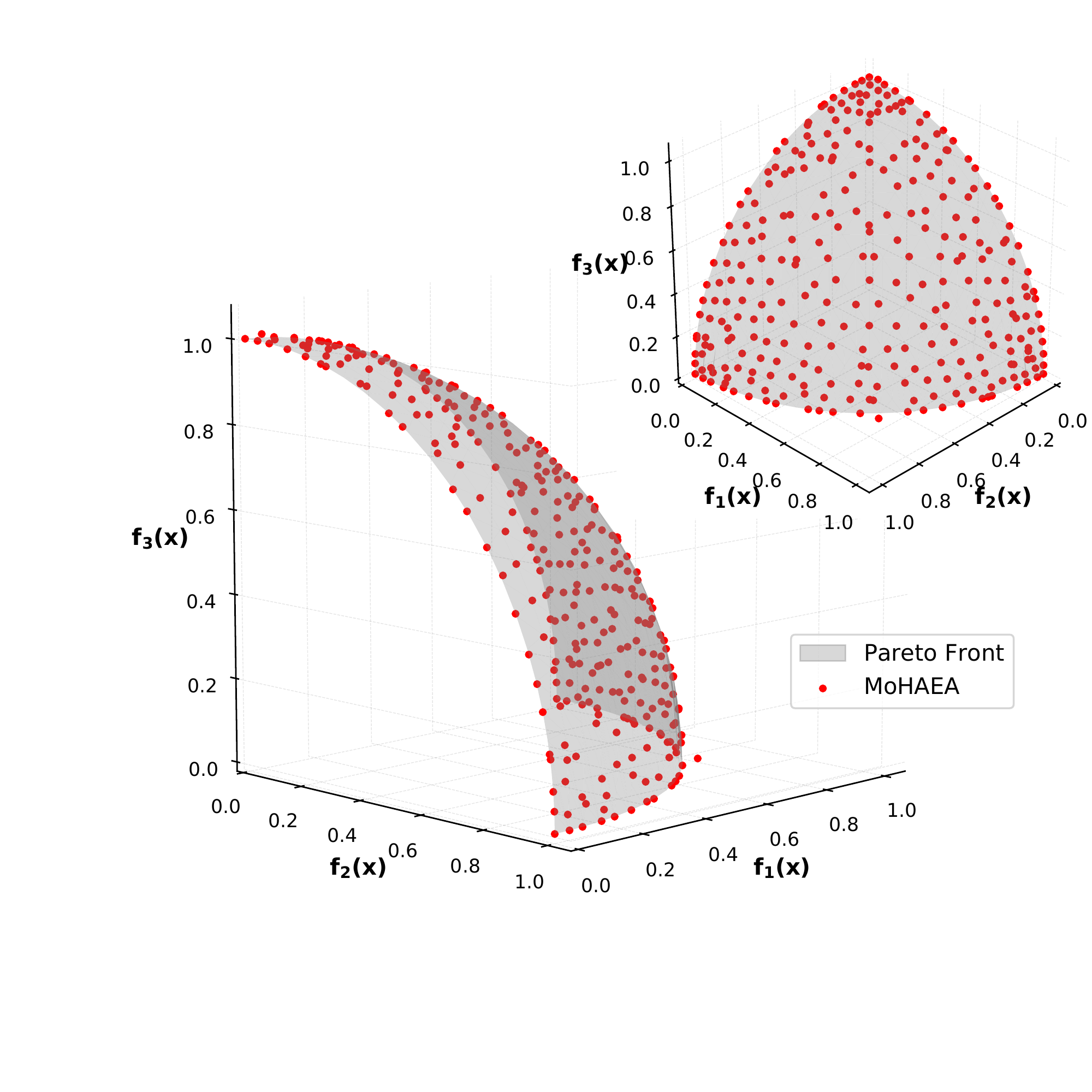}
        \caption{DTLZ4}
    \end{subfigure}
    \begin{subfigure}[t]{0.27\textwidth}
        \includegraphics[width=\textwidth, trim={0.5cm 2.5cm 0 0.5cm},clip]{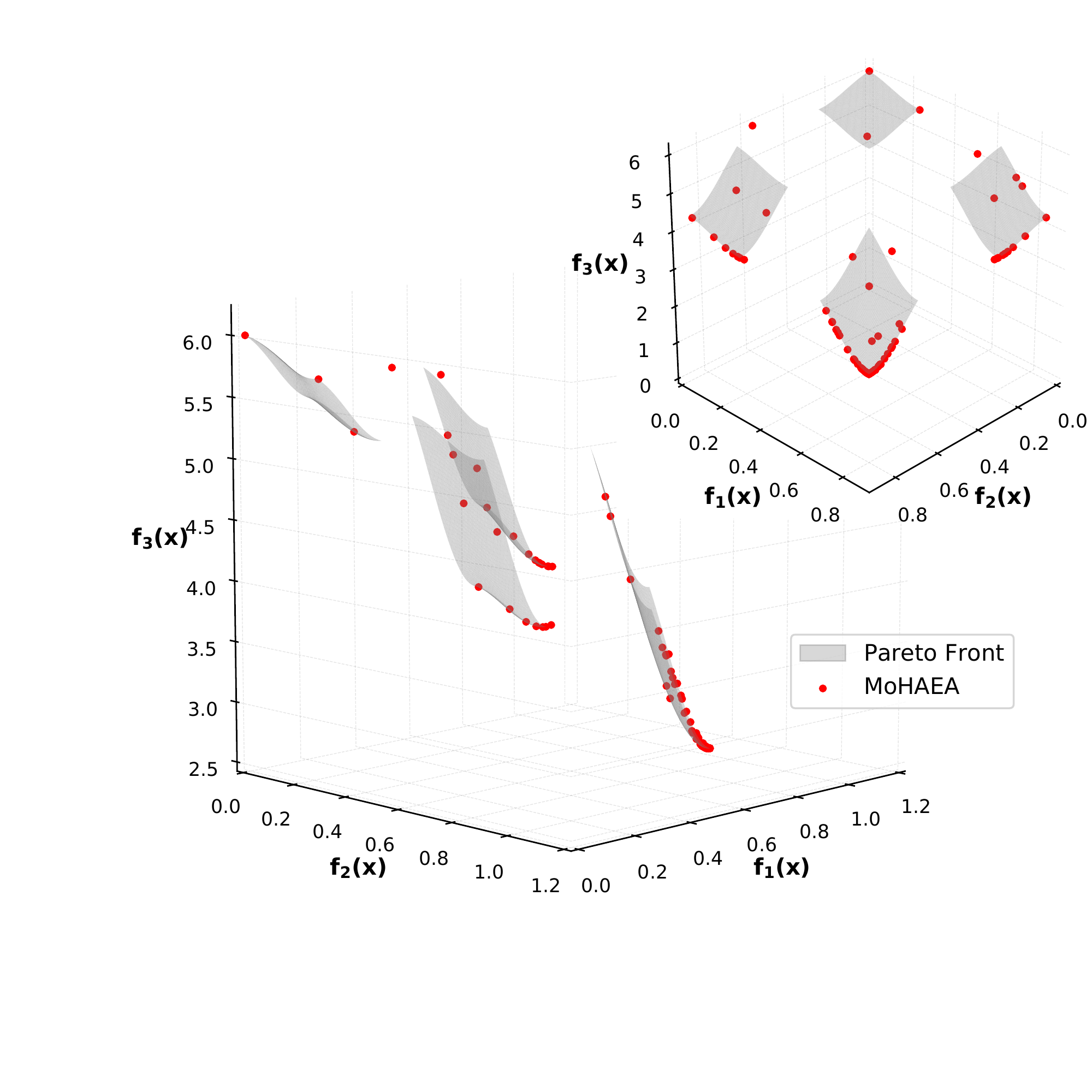}
        \caption{DTLZ6}
    \end{subfigure}
    \caption{Distribution of the non-dominated solutions with lowest IGD values found by MoHAEA on the ZDT and DTLZ problems.}
    \label{fig:ParetoFront}
\end{figure*}

\begin{figure*}[t]
    \centering
    \begin{subfigure}[t]{0.175\textwidth}
        \includegraphics[width=\textwidth, trim={0 0.4cm 0 0},clip]{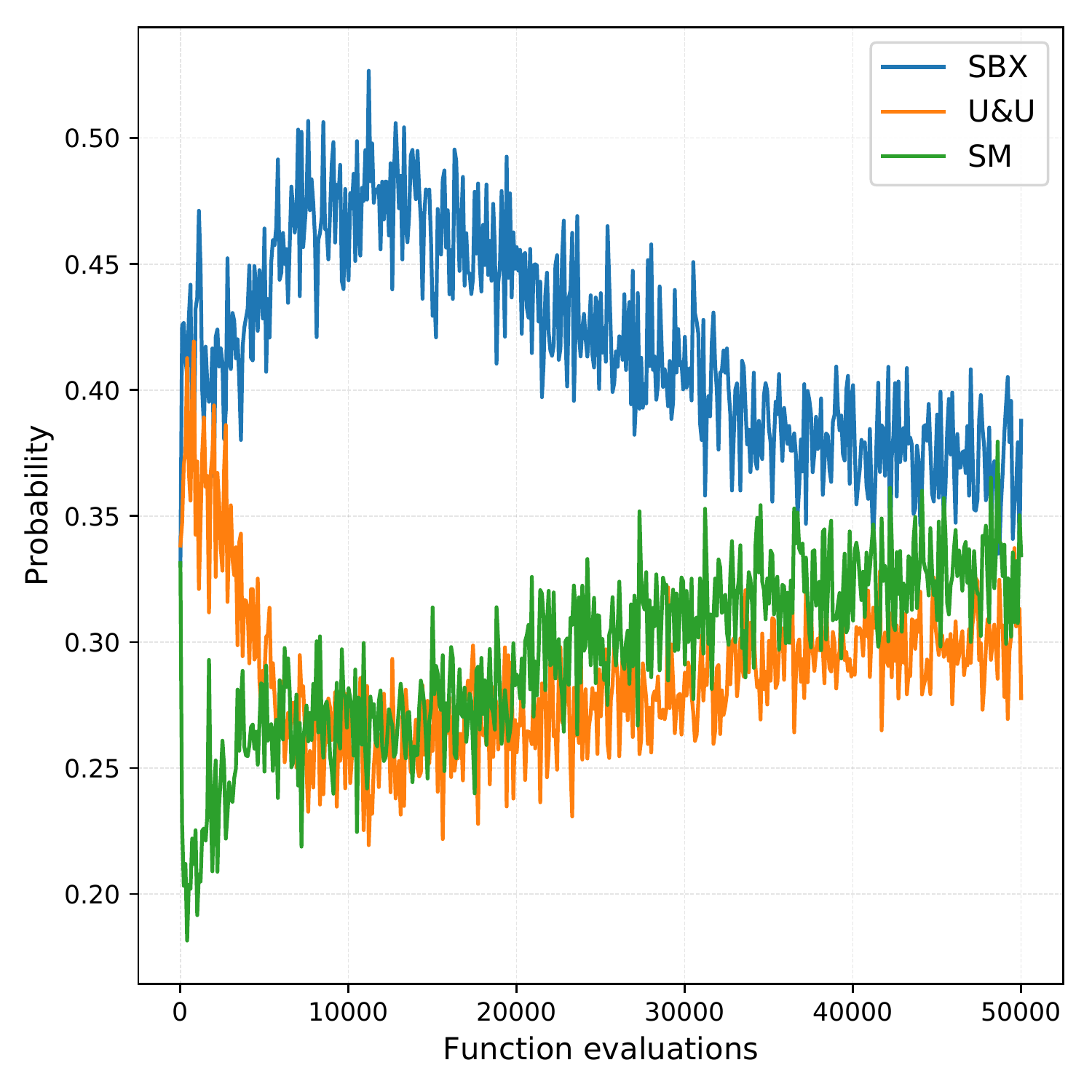}
        \caption{ZDT1}
    \end{subfigure}
    \hphantom{---}
    \begin{subfigure}[t]{0.175\textwidth}
        \includegraphics[width=\textwidth, trim={0 0.4cm 0 0},clip]{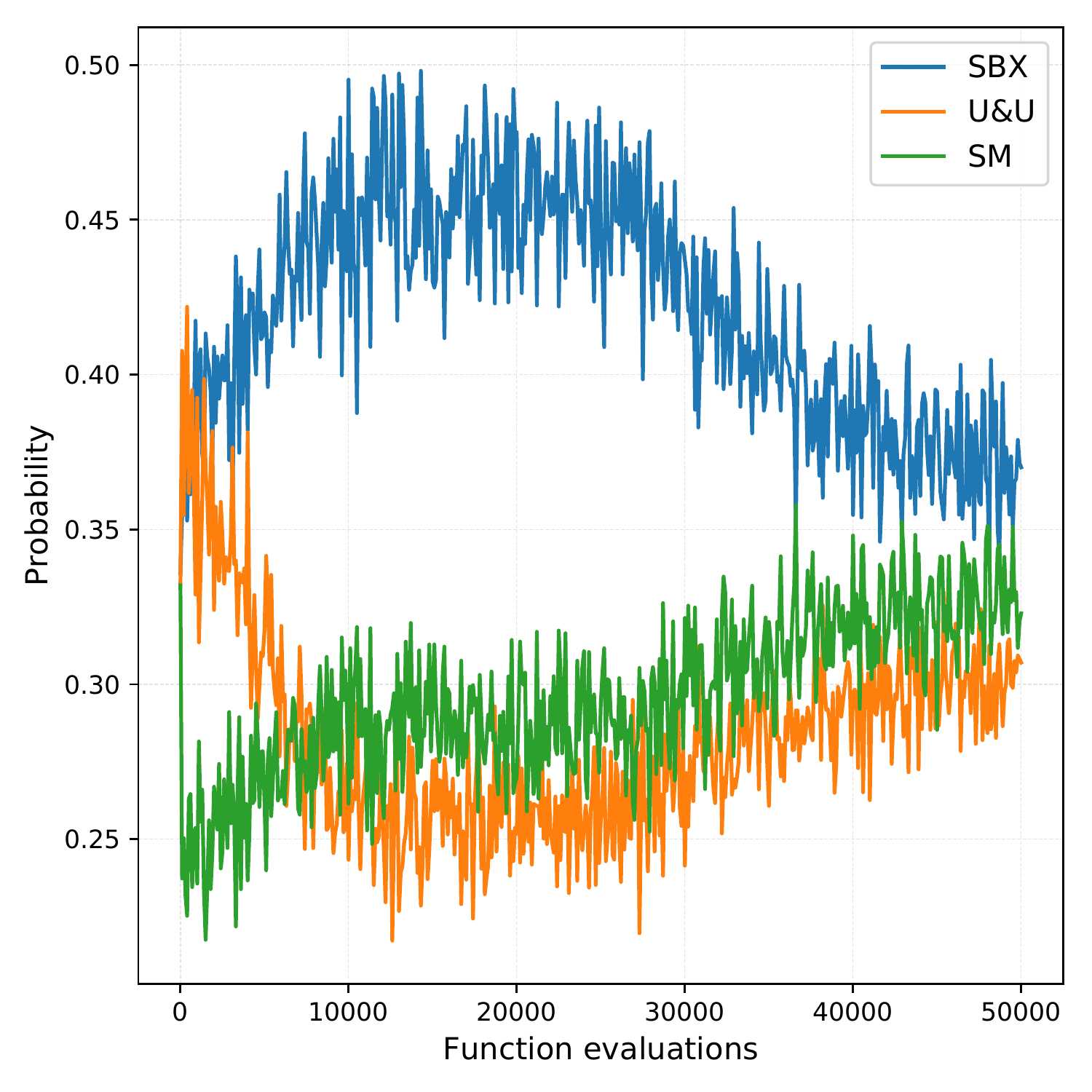}
        \caption{ZDT2}
    \end{subfigure}   
    \hphantom{---}
    \begin{subfigure}[t]{0.175\textwidth}
        \includegraphics[width=\textwidth, trim={0 0.4cm 0 0},clip]{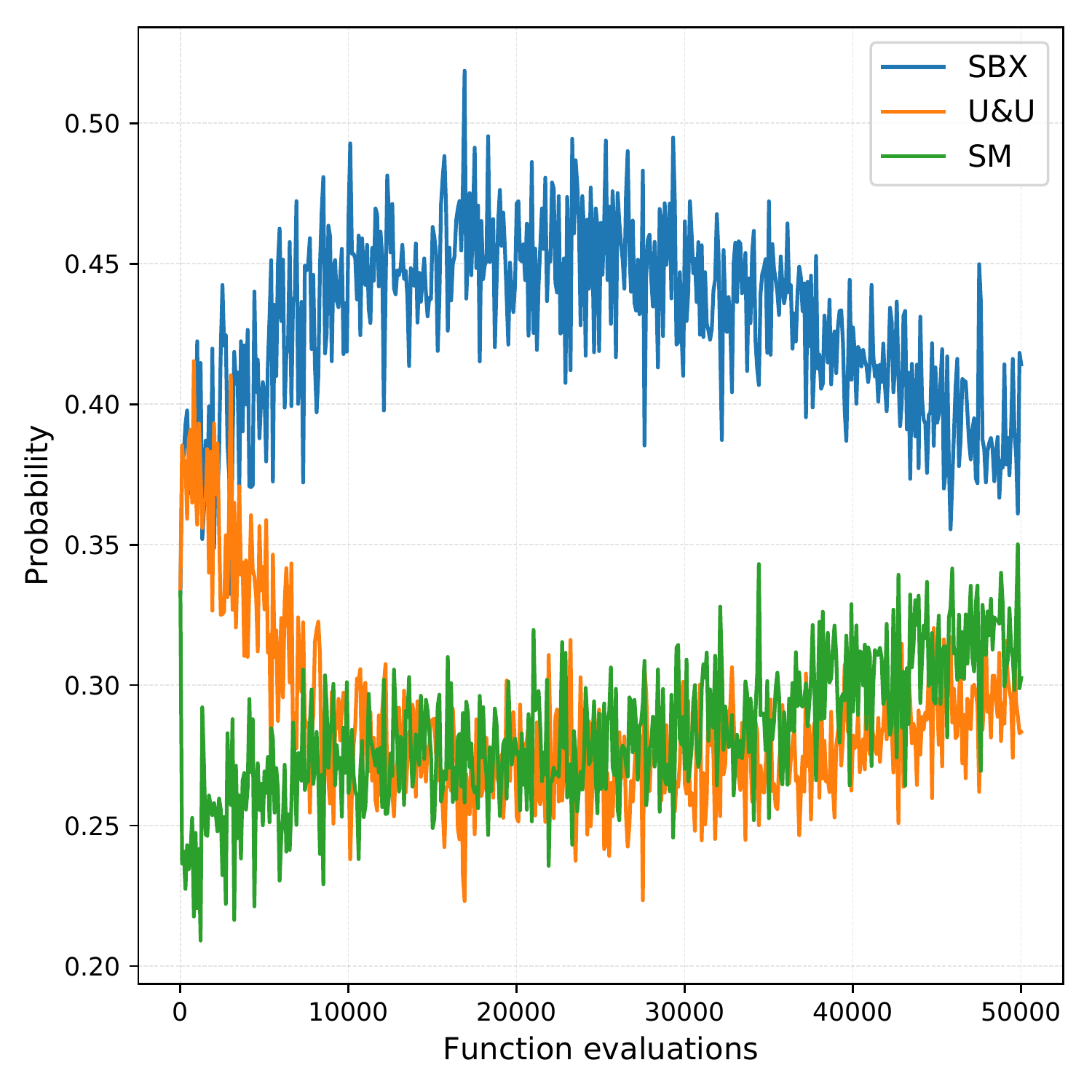}
        \caption{ZDT3}
    \end{subfigure}
    \hphantom{---}
    \begin{subfigure}[t]{0.175\textwidth}
        \includegraphics[width=\textwidth, trim={0 0.4cm 0 0},clip]{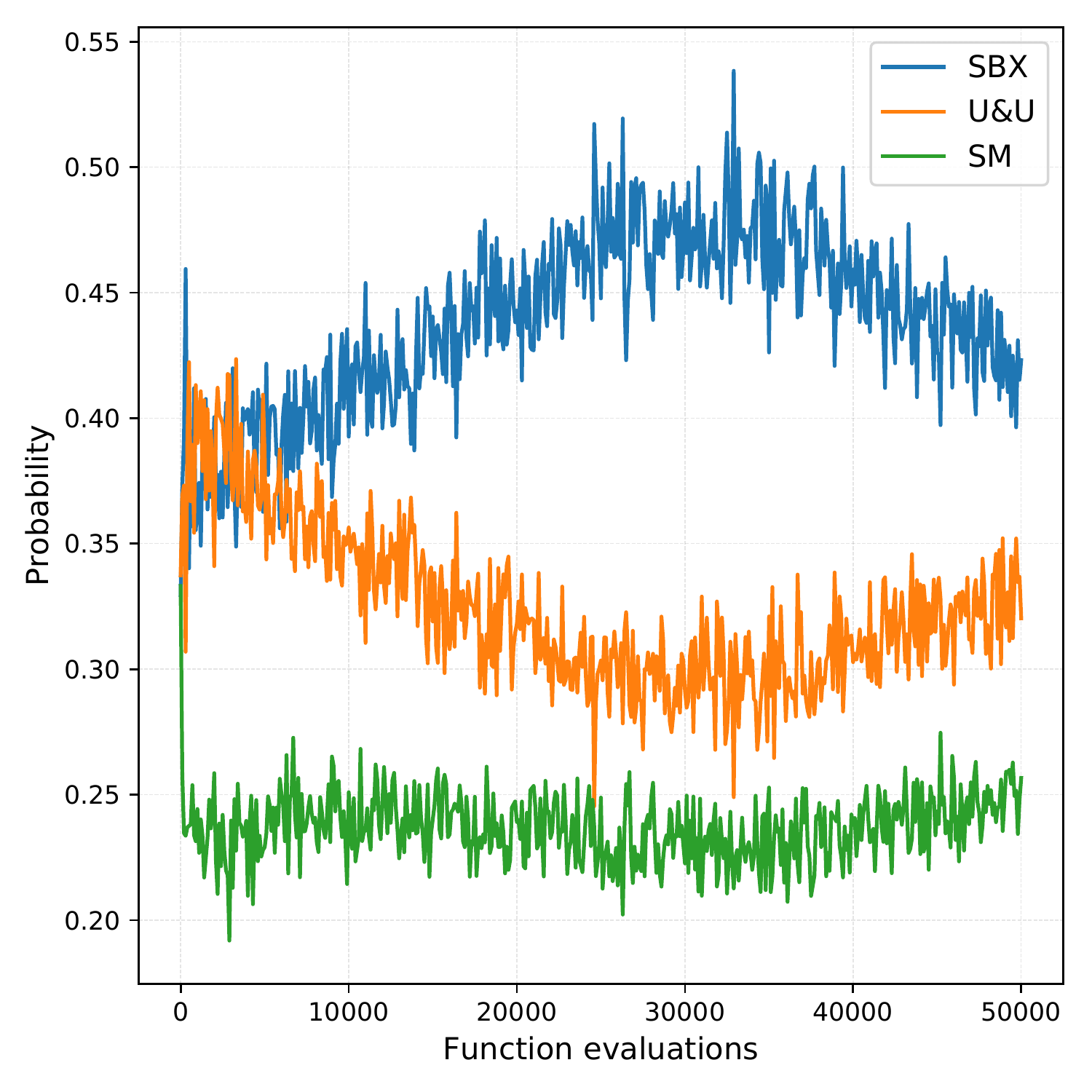}
        \caption{ZDT4}
    \end{subfigure}
    \begin{subfigure}[t]{0.175\textwidth}
        \includegraphics[width=\textwidth, trim={0 0.4cm 0 0},clip]{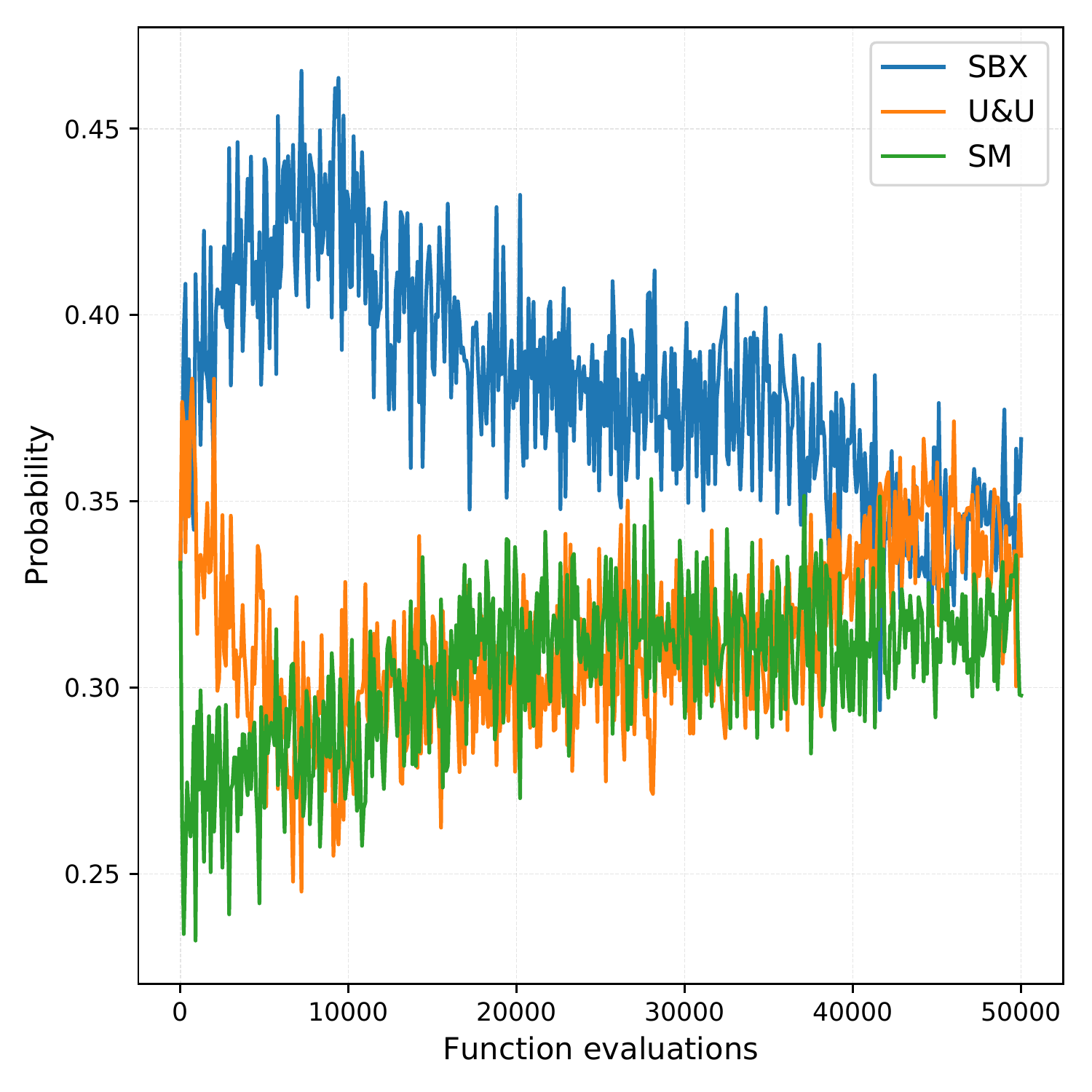}
        \caption{ZDT6}
    \end{subfigure}
    
    \begin{subfigure}[t]{0.175\textwidth}
        \includegraphics[width=\textwidth, trim={0 0.4cm 0 0},clip]{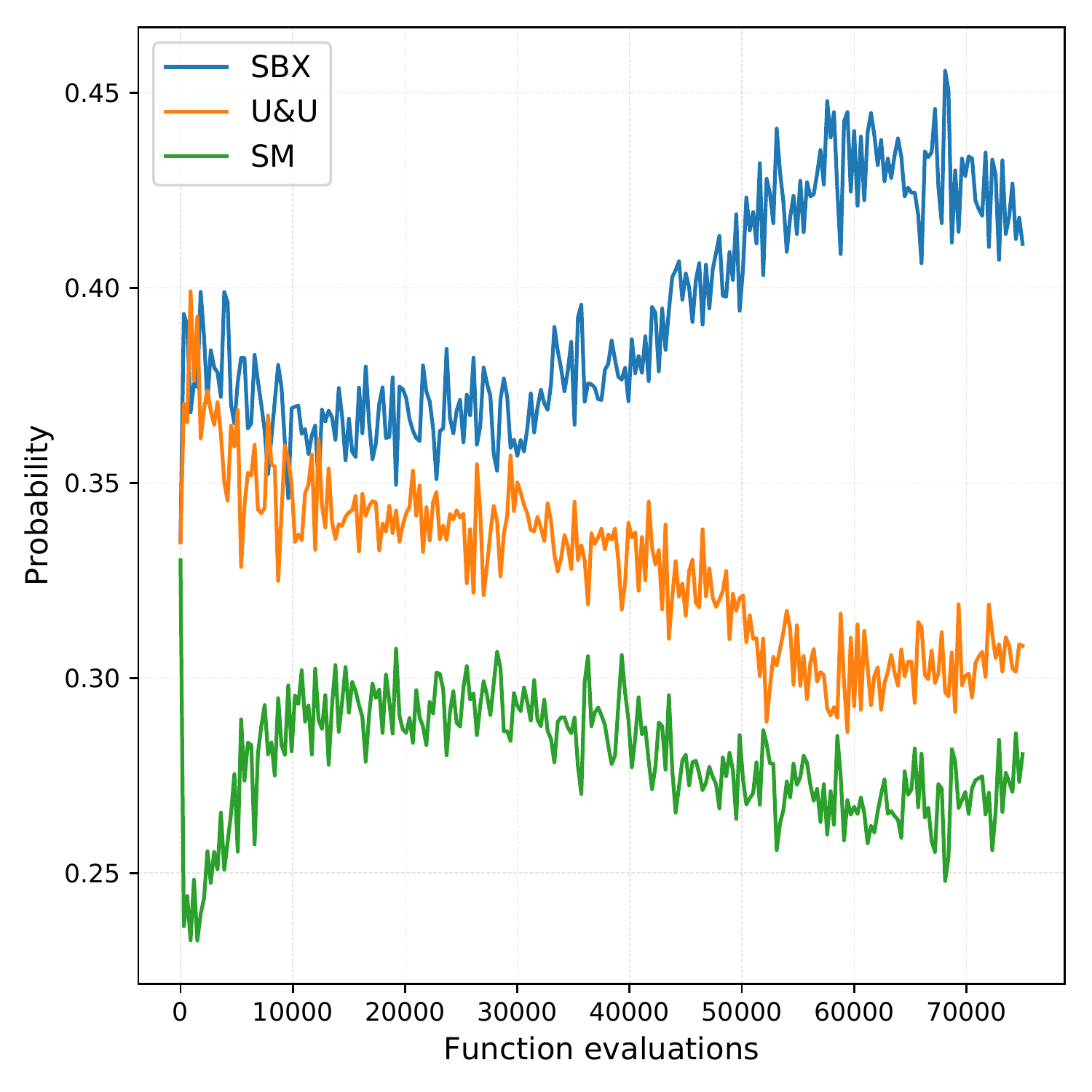}
        \caption{DTLZ1}
    \end{subfigure}
    \hphantom{---}
    \begin{subfigure}[t]{0.175\textwidth}
        \includegraphics[width=\textwidth, trim={0 0.4cm 0 0},clip]{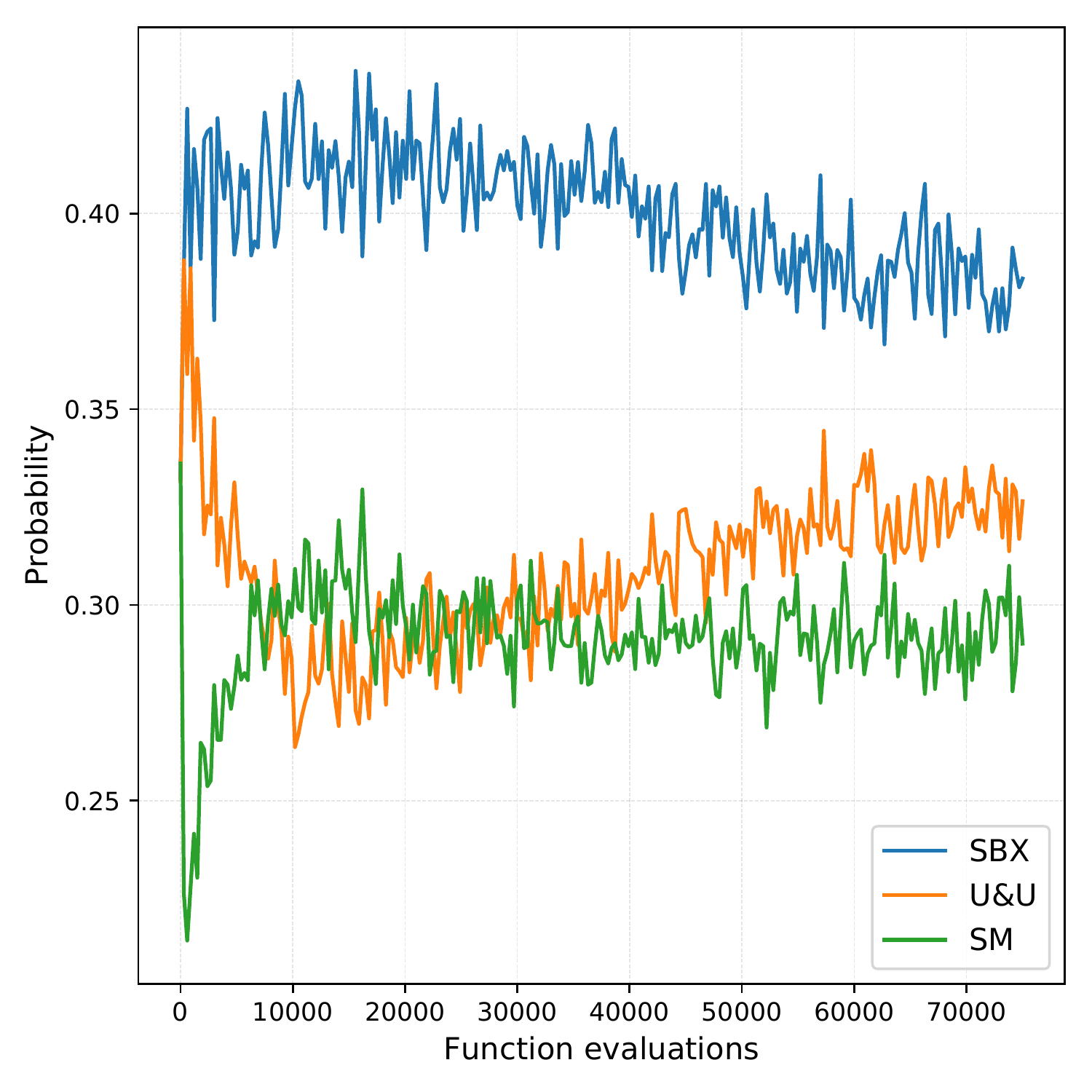}
        \caption{DTLZ2}
    \end{subfigure}   
    \hphantom{---}
    \begin{subfigure}[t]{0.175\textwidth}
        \includegraphics[width=\textwidth, trim={0 0.4cm 0 0},clip]{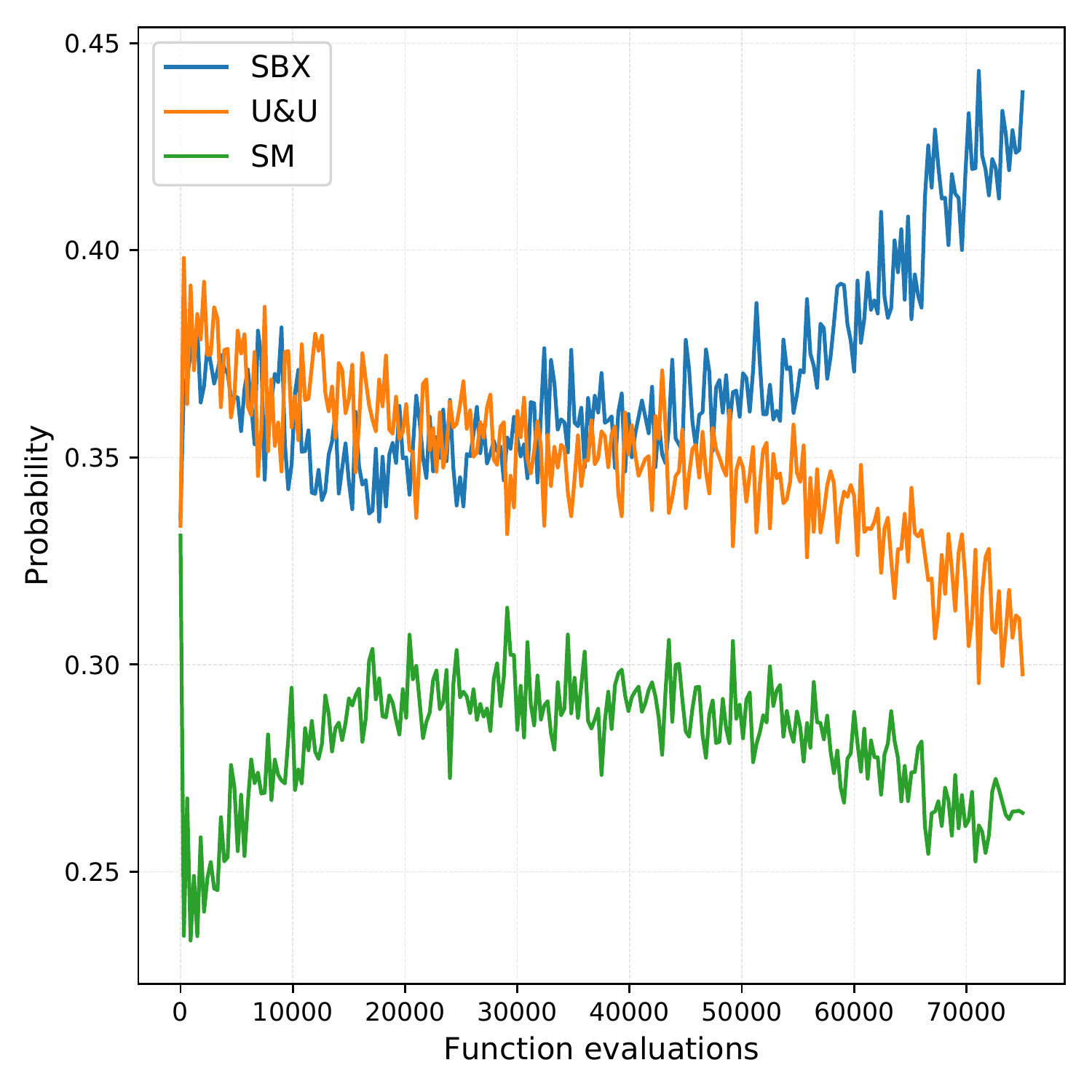}
        \caption{DTLZ3}
    \end{subfigure}
    \hphantom{---}
    \begin{subfigure}[t]{0.175\textwidth}
        \includegraphics[width=\textwidth, trim={0 0.4cm 0 0},clip]{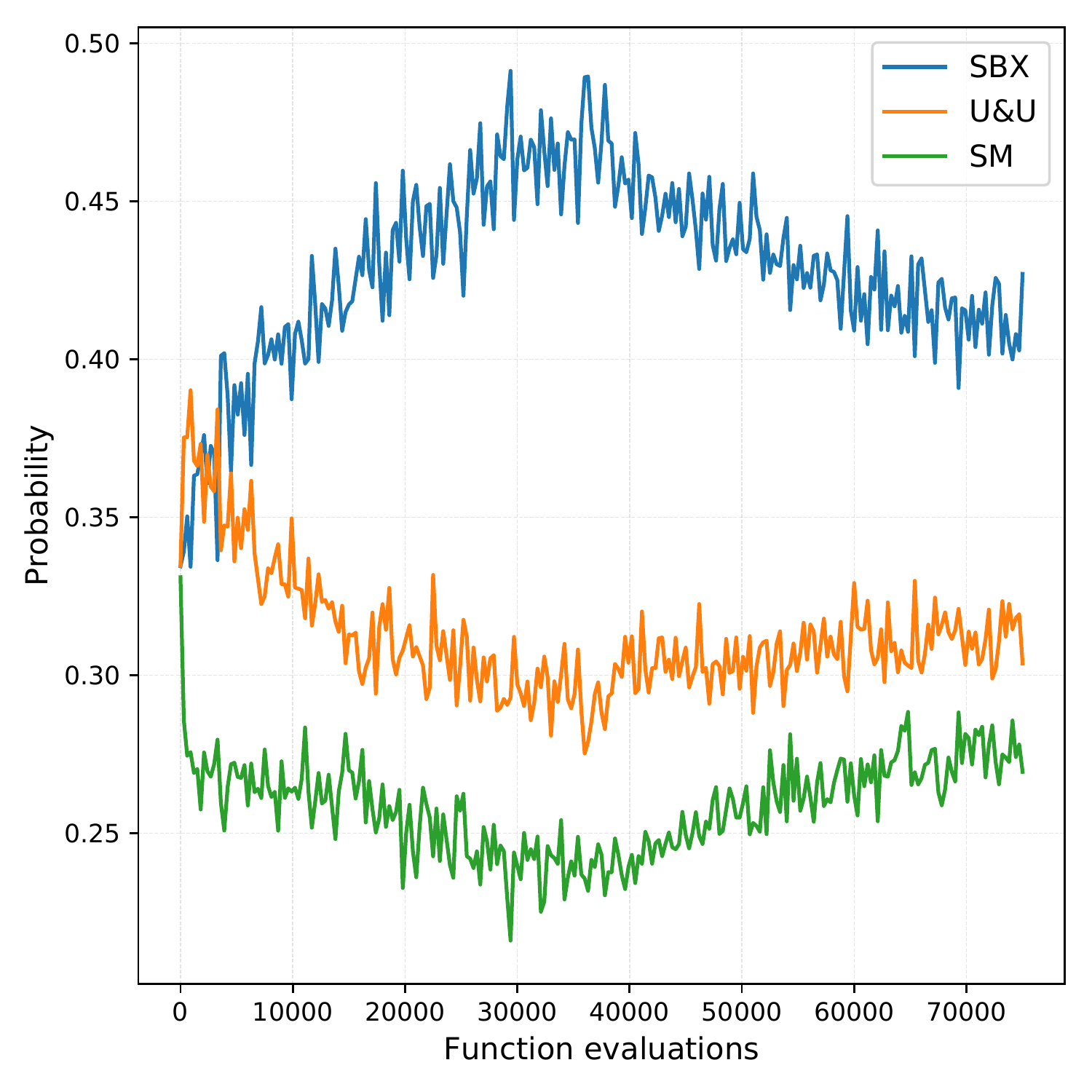}
        \caption{DTLZ4}
    \end{subfigure}
    \begin{subfigure}[t]{0.175\textwidth}
        \includegraphics[width=\textwidth, trim={0 0.4cm 0 0},clip]{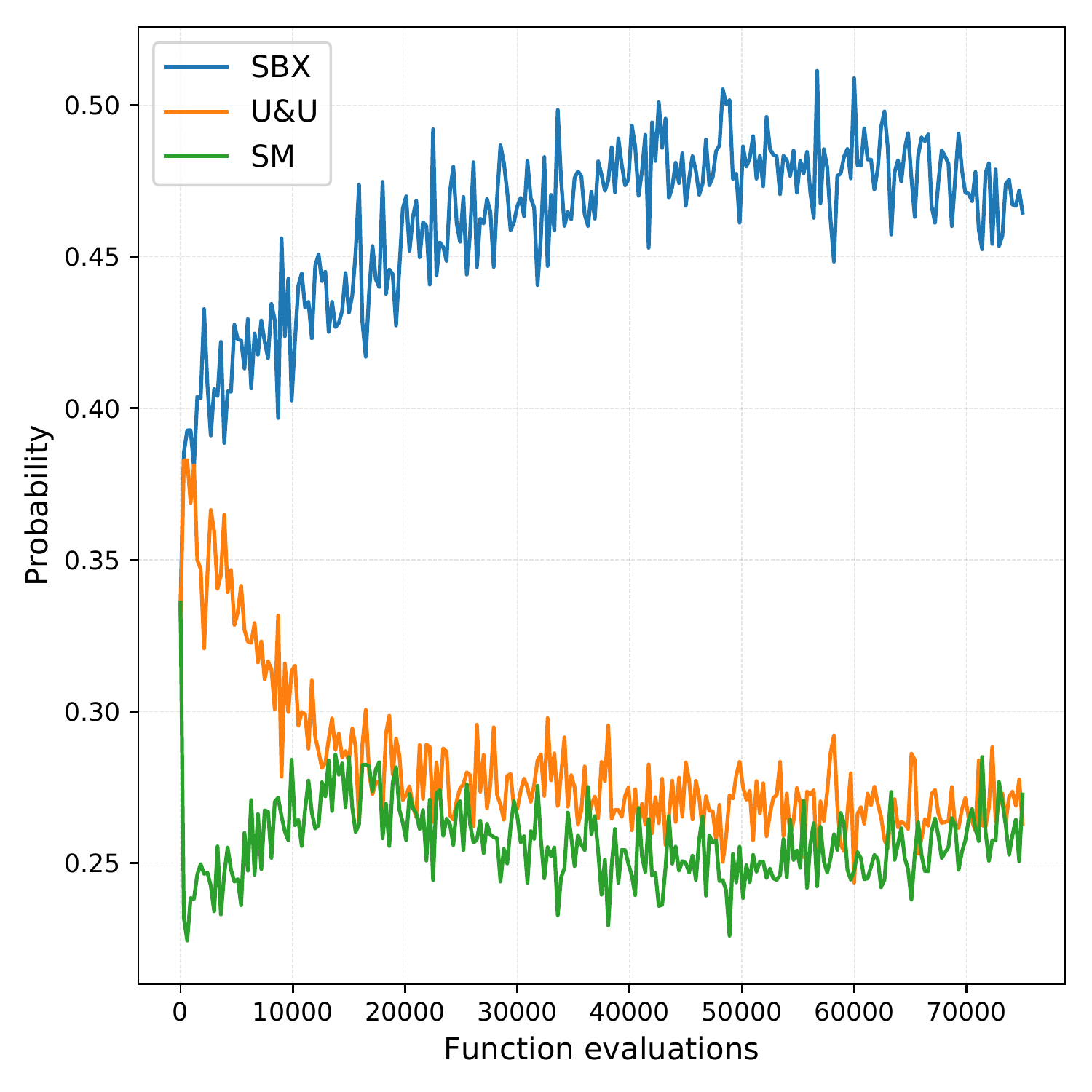}
        \caption{DTLZ6}
    \end{subfigure}
    \caption{Average of operators probabilities for  MoHAEA experiments over 30 different runs on the ZDT and DTLZ problems. Shrink Mutation SM (blue), Uniform \& Uniform U\&U (yellow), and Simulated Binary Crossover SBX (green).}
    \label{fig:OperatorRates}
\end{figure*}

\section{Results and Discussion}
\label{sec:Results}
\subsection{Experimental Results for ZDT and DTLZ Problems}    
Results obtained by MoHAEA are directly compared to those reported in \cite{qi2014moea} for MOEA/D, pa$\lambda$-MOEA/D, MOEA/D-AWA and NSGA-II, since the experimental setup is the same, including the way to compute IGD.

Table \ref{tb:comparison} presents the mean and computed deviation of the IGD metric values of the final solutions obtained by each algorithm on the ZDT problems and five DTLZ problems over 30 independent runs. This table reveals that in terms of the IGD metric, the final solutions obtained by MoHAEA are better than MOEA/D-AWA for 8 out of 10 problems, except for the ZDT3 and DTLZ3 problems. Especially, for the tri-objective DTLZ problems, MoHAEA (with PM or SM operator) performs much better than all of the compared algorithms. Moreover, MoHAEA(SM) outperforms the tested algorithms on the bi-objective ZDT problems, while MoHAEA(PM) works better on tri-objective DTLZ problems. This behavior can be explain as the high sensitivity of the combination of the genetic operators for different problems.

Figure \ref{fig:ParetoFront} shows, on the objective space, the distribution of the final solutions of MoHAEA(SM) with the lowest IGD value found by each algorithm for bi-objective ZDT and tri-objective DTLZ problems. Their mathematical descriptions and the ideal PFs of the ZDT and DTLZ problems can be found in \cite{zitzler2000comparison} and \cite{deb2002scalable}, respectively. It is visually clear that candidate solutions obtained by MoHAEA distribute significantly better, in terms of the uniformity of final solutions, than other tested algorithms \cite{qi2014moea}. 

\subsection{Learning Operator Rates}

Fig. \ref{fig:OperatorRates} shows the adaptation process of the genetic operator rates in MoHAEA (for MoHAEA(SM)) for the ZDT and DTLZ problems. As can be seen, Simulated Binary Crossover (SBX) and Uniform \& Uniform (U\&U) operator are more efficient than Shrink Mutation (SM) in the initial stages of the evolutionary. The rate performance shows that SBX is the best operator to be applied whenever the moment of the evolutionary process and for whatever problem (ZDT or DTLZ).

The property of MoHAEA allows showing when convergence is reached (i.e., ZDT1, ZDT2, and ZDT6 problems). So, no matter which operator is applied, the solution will not improve considerably.

\section{Conclusions and Future Work}
\label{sec:Conclusion}

In this paper, a multi-objective version of the Hybrid Adaptative Evolutionary Algorithm (HAEA) called (MoHAEA) is proposed. The proposed MoHAEA implements a unified paradigm, which combines dominance and decomposition-based approaches (Dominance penalty). Since each individual is evolved independently, the proposed dominance penalty approach defines a subproblem for fitness evaluation for each individual, reducing the extra parameters considerably for the algorithm (Neighborhood $T$ and penalty $\theta$).

The performance of MoHAEA has been investigated on a set of unconstrained well known ZDT and DTLZ instances, and compared with four other states of the art algorithms is made: MOEA/D, pa$\lambda$-MOEAD, MOEA/D-AWA, and NSGA-II. Based on the above comparisons, it can be concluded that the proposed MoHAEA can, in general, maintain a right balance between well-converged and well-distributed Pareto Front assisted by structured reference points. Indeed, through the test on problems with varying features, MoHAEA outperforms MOEA/D, pa$\lambda$-MOEA/D, MOEA/D-AWA, and NSGA-II on most of them in terms of IGD values. 

Future work includes the generalization of MoHAEA for multi-objective optimization with more than 3 objectives, and improvements on the reference points assignation at the initialization stage. Also, it will be interesting to investigate the performance of MoHAEA on a broader range of problems, such as problems with complicated PS shapes \cite{liu2013decomposition, saxena2011framework}, combinatorial optimization problems \cite{ehrgott2000survey}, and problems in real-world with a large number of objectives.

%


\bibliographystyle{ACM-Reference-Format}
\bibliography{sample-bibliography} 

\end{document}